\documentclass[12pt,a4paper,reqno]{amsart}
\usepackage{amsmath}
\usepackage{amsfonts}
\usepackage{amssymb}
\numberwithin{equation}{section}

\theoremstyle{plain}

\newtheorem{theorem}[subsection]{Theorem}
\newtheorem{proposition}[subsection]{Proposition}
\newtheorem{lemma}[subsection]{Lemma}
\newtheorem{corollary}[subsection]{Corollary}

\theoremstyle{definition}

\newtheorem{definition}[subsection]{Definition}

\newtheorem{question}[subsection]{Question}

\renewcommand{\leq}{\leqslant}
\renewcommand{\geq}{\geqslant}

\newsavebox{\proofbox}
\savebox{\proofbox}{\begin{picture}(7,7)%
  \put(0,0){\framebox(7,7){}}\end{picture}}

\newcommand{\md}[1]{\ensuremath{(\mbox{mod}\, #1)}}

\def\proof{\noindent\textit{Proof. }}
\def\endproof{\hfill{\usebox{\proofbox}}}

\def\scs{\scriptstyle}

\def\ni{\noindent}

\begin{document}

\title{Sum-free sets in abelian groups}

\author{Ben Green}
\address{Trinity College\\
     Cambridge\\
     CB2 1TQ\\
    England}
\email{bjg23@hermes.cam.ac.uk}

\author{Imre Z. Ruzsa}
\address{Alfr\'ed R\'enyi Mathematical Institute\\
     Hungarian Academy of Sciences\\
Budapest, Pf. 127\\
H-1364 Hungary}
\email{ruzsa@renyi.hu}

\thanks{The first author is supported by a Fellowship of Trinity College Cambridge.}

\begin{abstract}
\ni Let $A$ be a subset of an abelian group $G$ with $|G| = n$. We say that $A$ is sum-free if there do not exist $x,y,z \in A$ with $x + y = z$.
We determine, for any $G$, the maximal density $\mu(G)$ of a sum-free subset of $G$. This was previously known only for certain $G$.  We prove that the number of sum-free subsets of $G$ is $2^{(\mu(G) + o(1))n}$, which is tight up to the $o$-term. For certain groups, those with a small prime factor of the form $3k + 2$, we are able to give an asymptotic formula for the number of sum-free subsets of $G$. This extends a result of Lev, {\L}uczak and Schoen who found such a formula in the case $n$ even.
\end{abstract}

\maketitle

\section{Introduction and statement of results.} \ni Throughout this paper, $G$ will be a finite abelian group of order $n$. If $A$ is a subset of $G$ then we say that $A$ is \textit{sum-free} if there are no solutions to the equation $x + y = z$ with $x,y,z \in A$. Almost immediately upon making such a definition two natural questions present themselves:
\begin{question}\label{Q1} How big is the largest sum-free subset of $G$?
\end{question}
\begin{question}\label{Q2} How many sum-free subsets of $G$ are there?
\end{question}
We write $\mu(G)$ for the density of the largest sum-free subset of $G$, so that this subset has size $\mu(G)n$. We write $\mbox{SF}(G)$ for the set of all sum-free subsets of $G$. Observing that all subsets of a sum-free set are themselves sum-free, we have the obvious inequality
\[ |\mbox{SF}(G)| \; \geq \; 2^{\mu(G)n}.\]
Given this it is natural to introduce the notation 
\[ \sigma(G) \; = \; n^{-1}\log_2 |\mbox{SF}(G)|.\] Thus $\sigma(G) \geq \mu(G)$.\\[11pt]
A number of authors have addressed Questions \ref{Q1} and \ref{Q2}, and we take the opportunity to survey the best results known. \\[11pt]Interest in Question \ref{Q1} goes back over 30 years. Some straightforward observations get us quite a long way. Firstly note that $\mu(C_m) \geq \frac{1}{m}\left\lfloor \frac{m + 1}{3}\right\rfloor$, large sum-free sets being furnished by thinking of $C_m$ as $\mathbb{Z}/m\mathbb{Z}$ and taking appropriate intervals. It follows that $\mu(G) \geq 2/7$ whenever $G$ is cyclic. In fact the same inequality holds for all finite abelian groups, because $\mu(G) \geq \mu(G/H)$ for any quotient $G/H$ of $G$. Indeed if $\pi : G \rightarrow G/H$ is the canonical homomorphism and if $B \subseteq G/H$ is sum-free then so is the induced set $\pi^{-1}(B) \subseteq G$. This inequality is sharp, since Rhemtulla and Street \cite{RS} proved that $\mu(C_7^m) = 2/7$ for all $m$. \\[11pt]
By inducting from cyclic quotients in general one can easily prove the following. 
\begin{proposition}\label{P1}
Define a function $\nu$ from the set of all finite abelian groups to $[\frac{2}{7},\frac{1}{2}]$ as follows:
\begin{itemize}
\item If $n$ is divisible by a prime $p \equiv 2\,\mbox{\emph{(mod $3$)}}$ then $\nu(G) = \frac{1}{3} + \frac{1}{3p}$, where $p$ is the smallest such prime;
\item If $n$ is not divisible by any prime $p \equiv 2\,\mbox{\emph{(mod $3$)}}$, but $3|n$, then $\nu(G) = \frac{1}{3}$;
\item If $n$ is divisible only by primes $p \equiv 1\,\mbox{\emph{(mod $3$)}}$  then $\nu(G) = \frac{1}{3} - \frac{1}{3m}$, where $m$ is the exponent \emph{(}largest order of any element\emph{)} of $G$.\end{itemize}
Then $\mu(G) \geq \nu(G)$.
\end{proposition}
\ni It is convenient to have names for the three classes into which the finite abelian groups are divided by the above proposition.
\begin{definition}\label{typedef} If $n$ is divisible by a prime $p \equiv 2\md 3$ then we say that $G$ is type I. If $n$ is not divisible by any prime $p \equiv 2 \md 3$, but $3|n$, then we say that $G$ is type II. Otherwise, $G$ is said to be type III.
\end{definition}
\ni In the absence of obvious counterexamples it is natural to conjecture (cf. \cite{DY,Kedlaya}) that the lower bound of Proposition \ref{P1} is sharp, that is to say $\mu(G) = \nu(G)$. We prove this in the present paper.
\begin{theorem}\label{mainthm1}
We have $\mu(G) = \nu(G)$ for all finite abelian groups.
\end{theorem}
\ni This result has already been proved for type I and type II groups by Diananda and Yap \cite{DY}. It has also been proved for various type III groups, specifically groups of the form $C_{p^2} \times C_p$ and $C_{pq} \times C_p$ by Yap \cite{Yap1,Yap2} and elementary $p$-groups $C_p^m$ by Rhemtulla and Street \cite{RS}. Perhaps because of the increasing complexity of the proofs in these special cases, there does not appear to have been any progress on this problem since 1971.
We would recommend the paper of Kedlaya as an interesting introduction to this whole area of research.\\[11pt]
Moving on to Question \ref{Q2}, we remark that research in this direction was motivated by a conjecture of Cameron and Erd\H{o}s, now a theorem of the first author \cite{GreenCE}.
\begin{proposition} The number of sum-free subsets of $\{1,\dots,n\}$ is $O(2^{n/2})$.
\end{proposition} 
\ni An independent proof of this result was found somewhat later by A. Sapozhenko \cite{Sap2}.
Recently there has been progress on bounds for $|\mbox{SF}(G)|$ for various abelian groups. In the case $2|n$ an asymptotic was found by Lev, {\L}uczak and Schoen \cite{LLS} and independently by Sapozhenko \cite{Sap}:
\begin{proposition}[Lev-{\L}uczak-Schoen, Sapozhenko] \label{LLSthm} There is an absolute constant $c > 0$ such that \[ |\mbox{\emph{SF}}(G)| = (2^{\tau(G)} - 1)2^{n/2} + O(2^{(1/2 - c)n}),\]
where $\tau(G)$ is the number of even order components in the canonical decomposition of $G$ into a direct sum of cyclic groups.
\end{proposition}
\ni Lev, {\L}uczak and Schoen remark that they could take $c = 10^{-8}$, whilst Sapozhenko obtains the somewhat superior value $c = 0.017$. Even in the case $\tau(G)$ odd the bound of Proposition \ref{LLSthm}, whilst not giving an exact asymptotic, is decidedly non-trivial.\\[11pt]
The present authors \cite{GR}, improving on a result of Lev and Schoen \cite{LS}, found that when $p$ is a prime one has 
\[ \sigma(C_p) \; = \; 1/3 + o(1).\]
Given these results one might conjecture that $\sigma(G) = \mu(G) + o(1)$ for all $G$ (where $o(1)$ denotes a quantity that tends to zero as $n = |G| \rightarrow \infty$). This turns out to be the case, and is the second main result of our paper. 
\begin{theorem}\label{mainthm2} We have $\sigma(G)  =  \mu(G) + O\left((\log n)^{-1/45}\right)$. 
\end{theorem}
\ni For certain groups we are able to cast more light on the structure of a typical sum-free subset of $G$, and this leads to asymptotic bounds for $|\mbox{SF}(G)|$ rather than just for its logarithm. Let $p = 3k + 2$ be a prime. We say that $G$ is type I($p$) if it is type I and if $p$ is the \textit{least} prime factor of $n$ of the form $3k + 2$.
\begin{theorem}\label{mainthm3}
Suppose that $G$ is \textup{type I($p$)}. Then
\[ |\mbox{\emph{SF}}(G)| \; =  \; W \cdot \# \{\mbox{elements of order $p$}\} \cdot 2^{\mu(G)n}(1 + o_p(1)), \]
where $W = 1$ if $p = 2$ and $W = 1/2$ otherwise. 
\end{theorem}
\ni Observe that this does indeed generalise Proposition \ref{LLSthm}.

\section{An outline of the paper.}\label{sec2} \ni Ostensibly our paper contains two main strands: the determination of $\mu(G)$, the density of the largest sum-free subset of $G$, and the estimation of $\sigma(G)$. However, our strategy for counting sum-free sets means that these two strands are necessarily somewhat interlinked.\\[11pt]
Let us begin by saying a few words about this strategy. The main idea is to define, for any finite abelian group $G$, a certain family $\mathcal{F}$ of subsets of $G$. If $B \subseteq G$ is a set then we say that $(x,y,z) \in B^3$ is a \textit{Schur triple} if $x + y = z$.
\begin{proposition}\label{granprop} There is a family $\mathcal{F}$ of subsets of $G$ with the following properties:
\begin{enumerate}
\item $\log_2 |\mathcal{F}| = o(n)$;
\item Every $A \in \mbox{\emph{SF(}G\emph{)}}$ is contained in some $F \in \mathcal{F}$;
\item If $F \in \mathcal{F}$ then $F$ is \emph{almost sum-free}, meaning that $F$ has $o(n^2)$ Schur triples. 
\end{enumerate}
\end{proposition} \ni Proposition \ref{granprop} will be proved in \S 3. For the reader interested only in $\mu(G)$, this section can be completely ignored.\\[11pt]
In the later sections of the paper we will show that if $F$ is almost sum-free then $|F|$ cannot be much larger than $\mu(G)n$, the size of the largest sum-free set. This result is, perhaps, sufficiently important to be stated as a separate proposition.
\begin{proposition}\label{prop371} Suppose that a set $F \subseteq G$ has $\delta n^2$ Schur triples. Then 
\[ |F| \; \leq \; \left(\mu(G) + 2^{20}\delta^{1/5}\right)n.\] 
\end{proposition}
\ni This estimate, along with Proposition \ref{granprop}, immediately implies Theorem \ref{mainthm2}. Indeed, associate to each $A \in \mbox{SF}(G)$ some $F \in \mathcal{F}$ for which $A \subseteq F$. For a given $F$, the number of $A$ which can arise in this way is at most $2^{|F|}$. Thus we have the bound
\[ |\mbox{SF}(G)| \; \leq \; \sum_{F \in \mathcal{F}}2^{|F|} \; \leq \; |\mathcal{F}|\max_{F \in \mathcal{F}} 2^{|F|} \; = \; 2^{(\mu(G) + o(1))n}.\]
For groups of type I($p$) we will prove what amount to rough structure theorems for sets $F \in \mathcal{F}$ with size close to the maximal size $(\mu(G) + o(1))n$. This leads, by arguments similar to the above, to more precise counting results such as Theorem \ref{mainthm3}.\\[11pt]
We will prove Proposition \ref{prop371} using two slightly different arguments, one for groups of type I or II and the other for groups of type III.
It would be natural to try and show that if $F \subseteq G$ has $o(n^2)$ Schur triples then one may find a genuinely sum-free set $S \subseteq F$ with $|S| \geq |F| - o(n)$. Such a result is true, and is addressed in a preprint of the first-named author \cite{GreenReg}. However, the issues involved are rather complicated and the dependence between the $o(n^2)$ and the $o(n)$ coming from this approach seems to be extremely bad, certainly not as good as in Proposition \ref{prop371}.\\[11pt] 
When $|F| > n/3$ one can prove such an assertion by using an argument due to Lev, {\L}uczak and Schoen. We discuss this in \S 4, which is another part of the paper which the reader interested only in $\mu(G)$ may safely ignore.\\[11pt]
In \S \ref{sec7} we deal with groups of type III. Fortunately, our argument for determining $\mu(G)$ in this case is robust enough that it can be tweaked so as to cover almost sum-free sets as well. That is, we will show in a single argument that an almost sum-free subset of a type III group $G$ has cardinality at most $(\nu(G) + o(1))n$. Unfortunately, the method we use for showing that $\mu(G) = \nu(G)$ in this case is already rather unwieldly and the need to consider almost sum-free sets makes things look even more complicated.\\[11pt]
Let us conclude this section with some notation.  Write $\Gamma$ for the group of characters on $G$. If $f : G\rightarrow \mathbb{R}$ is a function then we define the Fourier transform of $f$ at $\gamma \in \Gamma$ by the formula
\[ \hat{f}(\gamma) \; = \; \sum_{x \in G} f(x) \gamma(x).\]
If $A \subseteq G$ then we will abuse notation by identifying $A$ with its characteristic function, allowing ourselves to talk, for example, of $\widehat{A}(\gamma)$.
If $A,B \subseteq G$ then we write $r(A,B,x)$ for the number of pairs $(a,b) \in A \times B$ with $a + b = x$. 
\section{Granular structure in groups.} \label{sec3} \ni This section contains a proof of Proposition \ref{granprop}, which is our main tool for counting sum-free sets. Let us remind ourselves of the statement of this proposition. In fact, the following is a rather more precise formulation than the earlier one:\\[11pt]
\textbf{Proposition $\mathbf{\ref{granprop}}^{\prime}$ }\emph{ Let $n$ be sufficiently large.  Then there is a family $\mathcal{F}$ of subsets of $G$ with the following properties:
\begin{enumerate}
\item $\log_2 |\mathcal{F}| \; \leq \; n(\log n)^{-1/18}$;
\item Every $A \in \mbox{\emph{SF}}(G)$ is contained in some $F \in \mathcal{F}$;
\item If $F \in \mathcal{F}$ then $F$ is \emph{almost sum-free}, meaning that $F$ has at most $n^2(\log n)^{-1/9}$ Schur triples.
\end{enumerate}}
\noindent Very roughly, the key idea will be to take a set $A \subseteq G$ and use it to construct a new set $A'$ which is much coarser that $A$ (being a collection of fairly large ``grains'') but which nonetheless contains fairly detailed information about the sumset $A + A := \{a_1 + a_2 : a_1,a_2 \in A\}$.\\[11pt]
There are two types of granular structure that we will consider. One is quite simple: if $L$ is an integer we say that a set is $L$-granular of coset type if it is a union of cosets of some subgroup $G_1 \leq G$ having size at least $L$. The other type of granularity is necessary for groups which do not possess many subgroups. Let $L$ be an integer and $d\in G$ be an element of order $m\geq L$. Partition $G$ as follows. We take each coset of the subgroup generated by $d$, split it into $\lfloor m/L\rfloor$ sets of type $x, x+d, ..., x+(L-1)d$ and one set of size less than $L$. There are many ways to do this, and we fix one of them for each $d$. A set which is the union of ``grains'' like this is called $L$-granular of progression type (note that the ``leftover'' sets of size less than $L$ are \textit{not} counted as grains).
\begin{lemma}[Granularization]\label{lem7.4.2}
 Let $A\subseteq G$ be a sum-free set and let $\epsilon \in (0,\frac{1}{2})$ be a real number. Let $L$ and $L'$ be positive integers satisfying
\begin{equation}\label{eq7.1} n \; > \; L' (10L/\epsilon)^{2^{34}\epsilon^{-8}}.\end{equation}
Then there is another set $A'$, \emph{(}the\emph{)} granularization of $A$, such that:
\begin{enumerate}
 \item $A'$ is either $L$-granular of progression type, or else $L'$-granular of coset type;\item
  $|A \setminus A'| \leq \epsilon n/4$; \item $A'$ has at most $\epsilon n^2/4$ Schur triples.\end{enumerate} 
 \end{lemma}  
\proof
We will define a certain set $P$, which in turn will be used to define $A'$. We will also consider the function $g : \Gamma \rightarrow [-1,1]$ defined by 
 \[ g(\gamma ) = \frac{1}{|P|} \sum_{b\in P} \gamma (b) . \] 
This is a normalised version of the Fourier transform $\widehat{P}(\gamma)$. The set $P$ will either be a subgroup $G_1$ of size $|G_1| \geq L'$, or will be of the form \[ P \; = \; \{ -(L-1)d, -(L-2)d, \dots, (L-1)d \}, \] 
where $d \in G$ and $\mbox{ord}(d) > 2L$.
These two cases will correspond to the two types of granular structure.
Taking $\delta = 2^{-16}\epsilon^4$, we will find a set $P$ of the above form so that  \begin{equation} \label{eq7.2}\left|\widehat{A}(\gamma) \left(1-g(\gamma)\right)  \right| \; \leq \;  \delta n \end{equation}    for all $\gamma \in \Gamma$. The function $g$, and particularly \eqref{eq7.2}, will be used to prove property (iii).\\[11pt]  
Let us now define $A'$.  If $P=G_1$, a subgroup, we let $A'$ be the union of those cosets of $P$ that contain at least $\epsilon |P |/4$ elements of $A$. Properties (i) and (ii) are clear in this case. If $P = \{ -(L-1)d, \dots , (L-1)d \}$ is an arithmetic progression, with difference $d$ having order $m$, then we consider the $L$-granular structure of progression type with common difference $d$. Let $A'$ be the union of those grains that contain at least $\epsilon L/8$ elements of $A$. Now $A \setminus A'$  contains at most $\epsilon L/8$ elements from each grain, making no more than $\epsilon n/8$ in total, plus at most $L$ elements from each of the $n/m$ ``leftover'' sets. Hence we have  \[  |A \setminus A' | \; \leq \; \epsilon n/8 + nL/m \; \leq \;  \epsilon n/4, \]   provided that  \begin{equation}\label{eq7.3}m \; \geq \;  8L/\epsilon  . \end{equation}    To establish property (iii) we consider an auxillary function $a_1$ defined by   
\[ a_1(x) \; = \; |A \cap (P + x) |/ |P|. \]  This is extremely natural since the Fourier transform of $a_1$ is just $\widehat{A}\cdot g$. Therefore we have, since $A$ is sum-free,
\begin{eqnarray}
\sum_{x + y = z} \nonumber a_1(x)a_1(y)a_1(z) & = & \sum_{x + y = z}\left(a_1(x)a_1(y)a_1(z) - A(x)A(y)A(z)\right)\\ & = & \nonumber
n^{-1} \sum_{\gamma}\left( |\widehat{a}_1(\gamma)|^2\widehat{a}_1(\gamma) - |\widehat{A}(\gamma)|^2\widehat{A}(\gamma)\right) \\ & = & \nonumber n^{-1}\sum_{\gamma} |\widehat{A}(\gamma)|^2\widehat{A}(\gamma)\left(1 - |g(\gamma)|^2g(\gamma)\right)\\ \nonumber  & \leq & n^{-1}\cdot\max_{\gamma} \left|\widehat{A}(\gamma)\right|\left|1 - |g(\gamma)|^2g(\gamma)\right|\cdot \sum_{\gamma}\left|\widehat{A}(\gamma)\right|^2 \\ \nonumber  & = & |A| \cdot \max_{\gamma} \left|\widehat{A}(\gamma)\right|\left|1 - g(\gamma)^3\right| \\
& \leq & \nonumber 3|A|\max_{\gamma} |\widehat{A}(\gamma)||1 - g(\gamma)|\\
 & \leq & 3\delta n^2,\label{eq71}
\end{eqnarray}
where the last two derivations use Parseval's identity and \eqref{eq7.2} respectively.
 Now consider an element $x \in A' $. If $P$ is a subgroup then $x+P$ contains at least $\epsilon |P |/4$ elements of $A$. When $P$ is a progression $P + x$ contains the grain of $A'$ to which $x$ belongs and hence at least $\epsilon|P|/16$ elements of $A$. In both cases $a_1(x)$ is at least $\epsilon/16$, and so $a_1(x) \geq \epsilon  A'(x)/16$ for all values of $x$. Thus, from \eqref{eq71}, we see that 
\begin{eqnarray} \nonumber
\mbox{number of Schur triples in $A'$} & = & \sum_{x + y = z} A'(x)A'(y)A'(z) \\ \label{eq72} & \leq & 2^{14}\epsilon^{-3}\delta n^2 \nonumber \\ & \leq & \epsilon n^2/4.\end{eqnarray} This completes the proof of property (iii).\\[11pt]
It remains to show that there is a set $P$ such that \eqref{eq7.2} holds. We also need \eqref{eq7.3} to hold if $P$ is a progression.
  Since $g(1)=1$ and $g(\gamma )\in [-1,1]$ for all $\gamma $, \eqref{eq7.2} automatically holds for $\gamma =1$ and for those $\gamma $ that satisfy $ |\widehat{A}(\gamma) |\leq \delta n/2$. Let $R $ be the set of all $\gamma  \ne  1 $ for which $ |\widehat{A}(\gamma) | > \delta  n/2$. We need to construct the set $P$ so that \eqref{eq7.2} holds for $\gamma \in R$.    Let $\Gamma _1$ be the subgroup generated by $R$ and let $G_1$ be the annihilator of $\Gamma _1$. If $ |G_1 |\geq L'$, we put $P=G_1$. In this case $g$ is a very simple function; indeed $g(\gamma )=1$ for $\gamma \in \Gamma_1$ (and it is 0 otherwise), so (\ref{eq7.2}) is immediate.    Assume then that $ |G_1 |<L'$. In order to find a suitable $d$ we reformulate condition (\ref{eq7.2}) in terms of the quantities $\arg \gamma (d)$.     Consider a general $\gamma \in \Gamma   $, and write $\arg \gamma (d) = \beta \in [-\pi ,\pi )$. We have 
\[ 1-g(\gamma ) \; = \; \frac{2}{2L-1} \sum_{j = 1} ^{L-1} (1 - \cos j\beta  ) \; \leq \; \frac{1}{2L-1} \sum_{j = 1}^{ L-1} ( j\beta  )^2 \; = \;\frac{L(L-1)}{6} \beta ^2 \; \leq \;  (L\beta )^2/6. \]
Hence a sufficient condition for \eqref{eq7.2} to hold is that 
\begin{equation}\label{eq7.5} |\arg \gamma (d) | \; \leq \; \frac{1}{L}\sqrt{\frac{6\delta n}{|\widehat{A}(\gamma) |}} \end{equation}    for all $\gamma \in R  $. We also need (\ref{eq7.3}) to hold. To achieve this we request that $d \notin G_1$ (so that $\gamma(d) \neq 1$ for at least one $\gamma \in R$) and 
strengthen condition (\ref{eq7.5}) to \begin{equation} \label{eq7.6}|\arg \gamma (d) |\;  \leq \; \frac{1}{L} \min \left ( \frac{\epsilon\pi}{4}, \sqrt {\frac{6\delta n}{|\widehat{A}(\gamma)|}} \right ) . \end{equation}  It follows by a standard application of the pigeonhole principle that we can find a $d \in G/G_1$ satisfying $d \neq 1$ and $ |\arg \gamma (d) | < \eta _\gamma $ for prescribed positive numbers $\eta _\gamma $ if  \begin{equation}\label{eq7.7} |G/G_1 |\; > \; \prod_{\gamma \in R}   \left( 1+ \lfloor 2\pi /\eta _\gamma \rfloor \right)  . \end{equation}     With $\eta _\gamma $ given by the right-hand side of equation (\ref{eq7.6}) we estimate the right-hand side of \eqref{eq7.7} from above as follows. It is at most 
\[  \prod_{\gamma \in R}  \left ( 1 + 8L \max \left ( \frac{1}{\epsilon},  \sqrt {\frac{ |\widehat{A}(\gamma) |}{\delta n}} \right ) \right ), \] which in turn is no more than   \[  (10L)^{|R|} \prod_{\gamma \in R}  \max \left ( \frac{1}{\epsilon}, \sqrt {\frac{ |\widehat{A}(\gamma)|}{\delta n}} \right ) . \]  To estimate this product we apply the following simple calculus lemma.
\begin{lemma}\label{lem7.3.2}
Let $x_1, ..., x_k$ be real numbers satisfying $x_i\ge 1$ and $\sum  x_i \leq K$, and assume $\tau \geq e^{1/e}$. Then we have \[\prod_{i=1}^k  \max (\tau , x_i) \le  \tau ^K .  \] \end{lemma}  \proof  For $x\ge 1$ and $\tau \geq e^{1/e}$ we have $\max (\tau , x) \leq \tau^x$ by calculus. \endproof\\[11pt]
Taking  $\tau = \epsilon^{-4}$ and $x_\gamma  = \left|\widehat{A}(\gamma)\right|^2\delta^{-2}n^{-2}$, we see using Parseval's identity that a suitable value of $K$ is $\delta^{-2}$. Thus we have \[ \prod_{\gamma \in R}  \max \left ( \frac{1}{\epsilon}, \sqrt { \frac{|\widehat{A}(\gamma) |}{\delta n}} \right ) \; \leq \;  \epsilon^{-\delta^{-2}}. \] Furthermore another application of Parseval's identity gives $|R| \leq 4\delta^{-2}$, and so the right-hand side of (\ref{eq7.7}) is at most $(10L/\epsilon)^{4\delta^{-2}}$. The left-hand side, however, is at least $n/L'$ and so the condition \eqref{eq7.1} implies that (\ref{eq7.7}) holds, and therefore that an element $d$ with the required properties can be found. This completes the proof of Lemma \ref{lem7.4.2}.\endproof\\[11pt]
It is a short step from Lemma \ref{lem7.4.2} to Proposition $\ref{granprop}^{\prime}$, which is the main result of this section. We prepare the ground with two very simple lemmas.
\begin{lemma}
\label{lem7.4.1} Suppose that $n$ is larger than some absolute constant and that $L \leq \sqrt{n}$. 
Then the number of subsets of $G$ which are $L$-granular \emph{(}of either coset or progression type\emph{)} is at most $2^{3n/L}$.
\end{lemma}
\proof The number of subgroups of $G$ is at most $2^{(\log_2 n)^2}$, since any subgroup may be generated by at most $\log_2 n$ elements. Thus the number of $L$-granular sets of coset type is at most $2^{(\log_2 n)^2 + n/L}$. Any $L$-granular set of progression type is associated with a partition of $G$ into at most $n/L$ grains, arising from the selection of an element $d \in G$ of order at least $L$. The number of such sets is thus at most $n2^{n/L}$. A short computation confirms the result.\endproof
\begin{lemma}\label{lem7.3.1} Suppose that $\rho$ is smaller than some absolute positive constant, and that $n$ is sufficiently large. Then the number of subsets of an $n$-element set of cardinality at most $\rho n$ is no more than $2^{n\sqrt{\rho}}$.
\end{lemma} 
\proof If $\rho < 1/n$ the result is trivial, so suppose this is not the case. The number in question is just $S = \sum_{k \leq \rho n} \binom{n}{k}$, which is certainly at most $n\binom{n}{\lfloor \rho n\rfloor}$ if $\rho < 1/2$. Using the well-known inequality $\binom{n}{k} \leq (en/k)^k$, we see that
\[ S \; \leq \; 2^{\rho\log_2(e/\rho)n + \log_2 n}.\]
Clearly $\rho\log_2(e/\rho) \leq \sqrt{\rho}/2$ for $\rho$ sufficiently small, and furthermore the fact that $\rho \geq 1/n$ guarantees that $\log_2 n \leq \frac{1}{2}n\sqrt{\rho}$. This completes the proof.\endproof\\[11pt]
Now set $L = L' = \lfloor \log n \rfloor$ and $\epsilon = (\log n)^{-1/9}$. One can easily check that, provided $n$ is sufficiently large, the condition \eqref{eq7.1} is satisfied. Thus we may apply Lemma \ref{lem7.4.2} with these values of $L,L'$ and $\epsilon$.\\[11pt]
Now for each $A \in \mbox{SF}(G)$ fix a set $A'$ (the existence of which is guaranteed by Lemma \ref{lem7.4.2}) and let $\mathcal{F}$ consist of all sets $A \cup A'$, for all $A \in \mbox{SF}(G)$. Then property (ii) of Proposition $\ref{granprop}^{\prime}$ is immediate from the construction of $\mathcal{F}$. To prove (iii), observe that $F = A \cup A'$ can be obtained from $A'$ by adding at most $\epsilon n/4$ elements of $A$. However the addition of a new $x \in G$ to some set $B$ cannot create more than $3n$ new Schur triples, and so any $F \in \mathcal{F}$ has at most $\epsilon n^2$ Schur triples. Finally, observe that by Lemmas \ref{lem7.4.1} and \ref{lem7.3.1} we have
\[ \log_2 |\mathcal{F}| \; \leq \; \frac{3n}{L} + \frac{n\sqrt{\epsilon}}{2} \; \leq \; n\sqrt{\epsilon} \; = \; n(\log n)^{-1/18}.\] This concludes the proof of Proposition $\ref{granprop}^{\prime}$.\endproof 
\section{A lemma of Lev, {\L}uczak and Schoen.} \ni In this section we apply a result of Lev, {\L}uczak and Schoen \cite{LLS} to show that if $F \subseteq G$ is suitably large and almost sum-free then $F$ has a large subset which is genuinely sum-free. The precise statment of this is given in Lemma \ref{LLSlem} below, but first we set up a piece of notation and recall a result from \cite{LLS}.
If $X \subseteq G$ we write $D = D_K(X)$ for the set of $K$-\textit{popular differences} of $X$, that is the set of all $d \in G$ which have at least $K$ different representations as $x_1 - x_2$, $x_1,x_2 \in X$. As usual, write $X - X := \{x_1 - x_2 : x_1,x_2 \in X\}$.
\begin{proposition}[\cite{LLS}, Proposition 1] \label{llsprop}With notation as above, suppose that $|D_K(X)| \leq 2|X| - 5\sqrt{K|X-X|}$. Then there is a subset $X' \subseteq X$ such that $|X \setminus X'| \leq \sqrt{K|X-X|}$ and $X' - X' \subseteq D_K(X)$.
\end{proposition}
\begin{lemma}\label{LLSlem} Let $\epsilon > 0$. Suppose that $F \subseteq G$ has cardinality at least $(1/3 + \epsilon)n$, and that $F$ has at most $\epsilon^3 n^2/27$ Schur triples. Then there is a set $S \subseteq F$ with cardinality $|S| \geq |F| - \epsilon n$ and which is sum-free.
\end{lemma}
\proof Set $N = \epsilon^3n^2/27$. If $N^2 \leq n$ then $N \leq \epsilon n$, and so we may simply remove an arbitrary element from each Schur triple to get our set $S$. Suppose then that $N^2 > n$, and set $K = \lceil N^{2/3}n^{-1/3}\rceil$. Observe that 
\begin{equation}\label{star5} N \; \geq \; |D_K(F) \cap F|K,\end{equation} since otherwise $F$ would have more than $N$ Schur triples. It follows that
\[ \frac{N}{K} \; \geq \; |D_K(F) \cap F| \; \geq \; |D_K(F)| + |F| - n.\] Therefore 
\[ |D_K(F)| \; \leq \; 2|F| - 3\epsilon n + \frac{N}{K} \; < \; 2|F| - 5\sqrt{Kn} \; \leq \; 2|F| - 5\sqrt{K|F - F|},\] the second inequality being a consequence of
\[ \frac{N}{K} + 5\sqrt{Kn} \; < \; 9(Nn)^{1/3} \; = \; 3\epsilon n.\]
By Proposition \ref{llsprop} there exists $F' \subseteq F$ such that $|F'| \geq |F| - \sqrt{Kn}$ and $F' - F' \subseteq D_K(F)$. The set $S = F' \setminus (F' - F')$ is clearly sum-free, and it is a subset of $F$. Moreover \eqref{star5} implies that
\[ |F' \cap (F' - F')| \; \leq \; |F \cap D_K(F)| \; \leq \; \frac{N}{K}\] and so we have the bound
\[ |S| \; \geq \; |F'| - \frac{N}{K} \; \geq \; |F| - \left(\frac{N}{K} + \sqrt{Kn}\right) \; > \; |F| - 3(Nn)^{1/3} \; = \; |F| - \epsilon n.\]
This concludes the proof of the lemma.\endproof\\[11pt]
The following corollary is a step in the direction of Proposition \ref{prop371}.
\begin{corollary}\label{cor44} Suppose that $F \subseteq G$ has $\delta n^2$ Schur triples. Then
\[ |F| \; \leq \; \left(\max\left(\frac{1}{3},\mu(G)\right) + 3\delta^{1/3}\right)n.\]
\end{corollary}
\proof Set $\epsilon = 3\delta^{1/3}$ in Lemma \ref{LLSlem}.\endproof

\section{Groups of type I and II.} \ni In this section we prove Theorems \ref{mainthm1} and \ref{mainthm2} for groups of types I and II (the reader may wish to recall Definition \ref{typedef}, in which these terms are defined), as well as Theorem \ref{mainthm3} (which applies to groups of type I($p$)). \\[11pt]
Theorem \ref{mainthm1} for groups of type I and II was established by Diananda and Yap, but we include a proof here for completeness (the reader may also wish to consult \cite{street-book}, which gives a comprehensive discussion of the state-of-the-art in 1972). 
A crucial ingredient is the following theorem of Kneser \cite{Kneser,Nath}:
\begin{proposition}[Kneser] Let $G$ be a finite abelian group, and suppose that $A$ and $B$ are subsets of $G$ with $|A + B| \leq |A| + |B| - X$ for some positive integer $X$. Then $A + B$ is the union of cosets of some subgroup $H \leq G$ with cardinality at least $X$.
\end{proposition}
\ni We apply this to sum-free sets through the following lemma.
\begin{lemma}\label{lem77}
Let $A \subseteq G$ be sum-free, let $r > 0$, and suppose that $|A| \geq n/3 + r$. Then there is a subgroup $H \leq G$, $|H| \geq 3r$, and a sum-free set $B \subseteq G/H$ such that $A$ is contained in $\pi^{-1}(B)$, where $\pi : G \rightarrow G/H$ is the canonical homomorphism.
\end{lemma}
\proof Since $A$ is sum-free, we have
\[ |A + A| \; \leq \; n - |A| \; \leq \; \frac{2n}{3} - r \; \leq \; 2|A| - 3r.\] Applying Kneser's theorem, we see that $A + A$ is a union of $H$-cosets, for some subgroup $H \leq G$ with $|H| \geq 3r$. Let $B = \pi(A)$.  To see that $B$ is sum-free, suppose that there are $x,y$ and $z \in B$ with $x + y = z$. By picking $h,h'$ so that $h + x$ and $h' + y$ both lie in $A$, one sees that $(A + A) \cap (H + z) \neq \emptyset$. However, $A + A$ is a union of $H$-cosets and therefore $A + A$ must contain \textit{all} of $H + z$. Thus it must intersect $A$, which is contrary to the assumption that $A$ is sum-free.\endproof
\begin{proposition}[Diananda--Yap]\label{prop496} We have $\mu(G) \leq \max(\nu(G),1/3)$. In particular, Theorem \ref{mainthm1} holds for groups of type I and II.\end{proposition}
\proof We use induction on $n = |G|$. Observe that the function $\nu$, as defined in Proposition \ref{P1}, has the property that $\nu(G) \geq \nu(G/H)$ for \textit{any} group $G$ and any $H \leq G$. Now suppose that $A \subseteq G$ is sum-free and that $|A| > \max(\nu(G),1/3)n$. Then $|A| > (n+1)/3$, so by Lemma \ref{lem77} there is a non-trivial subgroup $H \leq G$ and a sum-free set $B \subseteq G/H$ such that $A \subseteq \pi^{-1}(B)$. Hence, using the induction hypothesis, we have 
\[ |A| \leq \mu(G/H)n \leq \max( \nu(G/H),1/3)n \leq \max(\nu(G),1/3)n,\] a contradiction. Thus if $G$ is type I or II then $|A| \leq \nu(G)n$ which, by the constructions described prior to Proposition \ref{P1}, implies that $|A| = \nu(G)n$.\endproof\\[11pt]
Now that we know $\mu(G)$ for groups of type I and II, we can allow ourselves to move on to estimates for $|\mbox{SF}(G)|$. 
\begin{lemma}\label{prop371a} Proposition \ref{prop371} holds for groups of type I and II. That is, if $G$ is type I or II and if $A \subseteq G$ has $\delta n^2$ Schur triples then $|A| \leq \left(\mu(G) + 2^{20}\delta^{1/5}\right)n$.\end{lemma}
\proof For these groups we have $\mu(G) \geq 1/3$, and so the result (in fact, a rather stronger result) is immediate from Corollary \ref{cor44}.\endproof\\[11pt]
We have already sketched, in \S 2, how this result is relevant to Theorem \ref{mainthm2}.
\begin{proposition}\label{046} Theorem \ref{mainthm2} holds for groups of type I and II. That is, if $G$ is type I or II then we have $\sigma(G) = \mu(G) + O\left((\log n)^{-1/45}\right)$.
\end{proposition}
\proof Consider the family $\mathcal{F}$ of ``almost sum-free'' sets constructed in Proposition $\ref{granprop}^{\prime}$. Since every $F \in \mathcal{F}$ has no more than $n^2(\log n)^{-1/9}$ Schur triples we can infer from Lemma \ref{prop371a} that any such $F$ has cardinality at most \[ \mu(G)n + O(n(\log n)^{-1/45}).\] Recalling that $\log_2|\mathcal{F}| \leq n(\log n)^{-1/18}$ and that every $A \in \mbox{SF}(G)$ is contained in some $F \in \mathcal{F}$, the required estimate follows immediately.\endproof\\[11pt]
We will prove the same result for type III groups later on, and this will complete the proof of Theorem \ref{mainthm2}. For now, however, we complete our treatment of type I and II groups by obtaining Theorem \ref{mainthm3}, an asymptotic result for type I($p$) groups.\\[11pt]
The next lemma gives more detailed information about large sum-free subsets of type I($p$) groups. 
\begin{lemma}\label{prop29} Suppose that $G$ is type I\emph{(}$p$\emph{)} and write $p = 3k + 2$. Let $A \subseteq G$ be sum-free, and suppose that $|A| > \left(\frac{1}{3} + \frac{1}{3(p+1)}\right)n$. Then we may find a homomorphism $\psi : G \rightarrow \mathbb{Z}/p\mathbb{Z}$ such that $A$ is contained in $\psi^{-1}(\{k+1,\dots,2k+1\})$.
\end{lemma}
\proof By Lemma \ref{lem77} we know that there is a subgroup $H$, $|G/H| \leq p + 1$, such that $A$ is contained in $\pi^{-1}(B)$, where $\pi : G \rightarrow G/H$ is the canonical homomorphism and $B$ is a sum-free subset of $G/H$. We claim that $G/H \cong \mathbb{Z}/p\mathbb{Z}$. Indeed of all the possible quotients $G/H$ with $|G/H| \leq p + 1$, only the ones isomorphic to $\mathbb{Z}/p\mathbb{Z}$ are of type I. For all the others, $A$ would (by Proposition \ref{prop496}) have cardinality no more than $n/3$, which is contrary to assumption.\\[11pt]
We must also have $|B| = k + 1$. To classify the sum-free subsets of $\mathbb{Z}/p\mathbb{Z}$ having cardinality $k+1$, one can apply Vosper's theorem \cite{Nath,Vos}  detailing the cases for which equality holds in the Cauchy-Davenport inequality. Indeed we have
\[ |B + B| \; \leq \; p - |B| \; = \; 2|B| - 1,\]
and so Vosper's result guarantees that $B$ is an arithmetic progression of length $k+1$. It is easy to check that $B$, being sum-free, must be a dilate of $\{k+1,\dots,2k+1\}$.\endproof
\begin{lemma}\label{P7} Let $G$ be type I\emph{(}$p$\emph{)} and write $p = 3k + 2$.
With $o_p(2^{\mu(G)n})$ exceptions, all sum-free $A \subseteq G$ are described as follows. Take a homomorphism $\psi : G \rightarrow \mathbb{Z}/p\mathbb{Z}$, and take $A$ to be a subset of $\psi^{-1}(\{k+1,\dots,2k+1\})$ together with $o_p(n)$ further elements.
\end{lemma}
\proof Recall that $\mu(G) = \frac{1}{3} + \frac{1}{3p}$. Set $\delta = 1/6p^2$. Let us look once again at the family $\mathcal{F}$ constructed in Proposition $\ref{granprop}^{\prime}$. For each $F \in \mathcal{F}$, consider the collection of sum-free $A$ contained in $F$. The number of $A$ arising from any particular $F$ satisfying $|F| \leq (\mu(G) - \delta)n$ is $o_p(2^{\mu(G)}n)$, and hence so is the number of sets $A$ arising from \emph{all} such $F \in \mathcal{F}$. If $|F| \geq (\mu(G) - \delta)n$ then we may, by Lemma \ref{LLSlem}, find a sum-free set $S \subseteq F$ with $|S| = |F| - o_p(n)$. Lemma \ref{prop29} then tells us that $S$ is contained in $\psi^{-1}(\{k+1,\dots,2k + 1\})$ for some suitable homomorphism $\psi : G \rightarrow \mathbb{Z}/p\mathbb{Z}$. The lemma follows immediately.\endproof
\begin{lemma}\label{P7a} Let $G$ be type I\emph{(}$p$\emph{)} with $p = 3k + 2$.
With $o_p(2^{\mu(G)n})$ exceptions, all sum-free $A \subseteq G$ are described as follows. Take a homomorphism $\psi : G \rightarrow \mathbb{Z}/p\mathbb{Z}$, and take $A$ to be a subset of $\psi^{-1}(\{k+1,\dots,2k+1\})$.
\end{lemma}
\proof For brevity write $M_k = \{k+1,\dots,2k + 1\}$, considered as a subset of $\mathbb{Z}/p\mathbb{Z}$. By the last lemma, it suffices to look at sum-free sets which are ``almost'' induced from $\mathbb{Z}/p\mathbb{Z}$, being a subset of some $\psi^{-1}M_k$ together with $o_p(n)$ elements. Let $H = \ker \psi$, so that with the exception of $o_p(n)$ elements $A$ is contained in $(H + k+1) \cup \dots \cup (H + 2k + 1)$. Fix $i \notin M_k$ and $x \in H + i$ and count the number of $A$ containing $x$.
Observe that any such $i$ can be written as $y \pm z$ for some distinct $y,z \in M_k$. Split $(H +y) \cup (H + z)$ into $|H|$ pairs $(s,t)$ with $s \pm t = x$. In each of these pairs at most one element can lie in $A$, giving at most $3^{|H|}$ possibilities for $A \cap ((H + y) \cup (H + z))$. 
There are at most $2^{(k-1)|H|}$ possibilities for $A \cap \bigcup_{j \in M_k \setminus \{y,z\}} (H + j)$ and $2^{o_p(n)}$ choices for $A  \setminus \bigcup_{j \in M_k} (H + j)$. This means that
the number of $A \in \mbox{SF}(G)$ containing $x$ does not exceed 
\[ 3^{|H|}2^{(k-1)|H| +  o_p(n)} \; = \; 2^{n(k - 1 + \log_2 3 + o_p(1))/p}.\]
However, since $\mu(G) = (k+1)/p$, we have
\[ \frac{n}{p}(k - 1 + \log_2 3 + o_p(1)) \; < \; (\mu(G) - \textstyle \frac{1}{10p})n.\]
Summing over $i$ and $x$ confirms the lemma.\endproof\\[11pt]
\noindent\textit{Proof of Theorem \ref{mainthm3}.}  We now know that most sum-free sets are contained in $\psi^{-1}(M_k)$ for some non-zero homomorphism $\psi : G \rightarrow \mathbb{Z}/p\mathbb{Z}$. We need to understand sum-free sets $A$ which arise in this way from two distinct homomorphisms $\psi_1$ and $\psi_2$. Observe that $\psi_i^{-1}(M_k)$ is a union of $k+1$ cosets of $\ker \psi_i$. If $\ker \psi_1 \neq \ker \psi_2$ then $\psi_1^{-1}(M_k) \cap \psi_2^{-1}(M_k)$ has cardinality $(k+1)^2n/p^2$. The number of $A$ which are contained in such an object is certainly $o_p(2^{\mu(G)n})$, since the number of choices for $\psi_1$ and $\psi_2$ is at most $n^2 = 2^{o_p(n)}$. If $\ker \psi_1 = \ker \psi_2$ then either $\psi_1^{-1}(M_k) = \psi_2^{-1}(M_k)$ or $|\psi_1^{-1}(M_k) \cap \psi_2^{-1}(M_k)| \leq kn/p$. In the latter case there are again $o_p(2^{\mu(G)n})$ possibilities for $A$. In the former case $\psi_1\psi_2^{-1}$ defines an isomorphism from $\mathbb{Z}/p\mathbb{Z}$ to itself which preserves $M_k$. If $p = 2$ such an isomorphism must be the identity, but if $p \geq 5$ there are two such maps, the other being $x \mapsto -x$.\\[11pt]
Thus almost all $A \in \mbox{SF}(G)$ arise from just one equivalence class of homomorphisms $\psi$ in the manner detailed by Lemma \ref{P7a}, where we say that $\psi_1$ and $\psi_2$ are equivalent if either $\psi_1 = \psi_2$ or $\psi_1 = - \psi_2$. It follows that
\[ |\mbox{SF}(G)| \; = \; W \cdot \# \{\mbox{non-zero homomorphisms $\psi : G \rightarrow \mathbb{Z}/p\mathbb{Z}$}\} \cdot 2^{\mu(G)n}(1 + o_p(1)),\] where $W = 1$ if $p = 2$ and $W = 1/2$ if $p \geq 5$. This completes the proof of Theorem \ref{mainthm3}. Indeed, the number of non-trivial homomorphisms from $G$ to $\mathbb{Z}/p\mathbb{Z}$ equals the number of non-trivial homomorphisms from $\Gamma$ to the multiplicative group of complex $p$th roots of unity, which is the number of elements of $G$ of order $p$.\endproof
\section{A Pollard-Kneser result.}\label{sec6} \ni In this section we prove a result, Proposition \ref{KP}, which is necessary for parts of \S \ref{sec7}. The reader who is only interested in $\mu(G)$, and not in counting sum-free sets, may ignore this section. The result generalises (a weak version of) Kneser's theorem and a theorem of Pollard \cite{Pol}. It is of interest in its own right, and furthermore is a fairly simple modification of well-known results. Nonetheless, we have not been able to find a reference for it.\\[11pt]
Given a finite abelian group $G$ define $D(G)$ to be the size of the largest proper subgroup of $G$. We will often refer to this as the \textit{defect} of $G$.
\begin{proposition}\label{KP} Write $D = D(G)$. Suppose that $A$ and $B$ are subsets of $G$ with cardinalities $k$ and $l$ respectively, and suppose that $t \leq \min(k,l)$ is a non-negative integer. Then
 \[ \sum_x \min\left(t,r(A,B,x)\right) \;\geq \; t\min(n, k + l  - D - t).\]
 \end{proposition}
 \proof Suppose without loss of generality that $l\leq k$. Writing $S(l)$ for the statement that the theorem holds for a given value of $l$ and all $k \geq l$, we proceed by induction on $l$. To check the base case $S(1)$ observe that the theorem is always true when $t = l$. Indeed we always have $r(A,B,x) \leq l$, and so
 \[
 \sum_x \min\left(l,r(A,B,x)\right) \; = \; \sum_{x} r(A,B,x) \; = \; kl.\]
 Suppose then that $l \geq 2$, and that $S(l')$ whenever $l' < l$. If $k + l - D - t > n$ then, by the pigeonhole principle, we have $r(A,B,x) \geq D + t \geq t$ for all $x$. Thus $
\sum_x \min\left(t,r(A,B,x)\right) =  tn$,
so that the result is true in this case too. Suppose then that $k + l - D - t \leq n$. By translating $B$ if necessary we may assume that $0 \in B$. We now distinguish two further cases.\\[11pt]
Case 1: $A + B = A$. Then $A$ is a union of cosets of some subgroup $H \leq G$ which contains $B$, and so $D \geq l$. Furthermore each element of $A + B$ is represented exactly $l$ times: therefore $r(A,B,x) \geq t$ for all $x \in A + B$, and it follows that
\[
\sum_x \min(t,r(A,B,x)) \; = \; |A + B|t \; = \; kt \;\geq \; t(k + l - D - t),
\]
confirming the theorem in this case. Observe that this case is the only important difference between general groups and groups of prime order.\\[11pt]
Case 2: $A + B \neq A$. Choose elements $a \in A,b \in B$ with $a + b \notin A$. Replacing $A$ with $A - a$ we may now assume that $0 \in A$ and that $B \not\subseteq A$. Write $U = A \cup B$ and $I = A \cap B$, and observe that $0 < |I| < l$. Finally set $A' = A \setminus I$ and $B' = B \setminus I$. We have
\begin{equation}\label{star1} r(A,B,x) \; = \; r(U,I,x) + r(A',B',x).\end{equation}
If $1 \leq t \leq |I|$ then the theorem follows immediately, since by (\ref{star1}) and the induction hypothesis
\begin{eqnarray*}
\sum_x \min(t,r(A,B,x)) & \geq & \sum_x \min(t,r(U,I,x)) \\ & \geq & t\min(n, |A| + |B | - D - t).
\end{eqnarray*}
In the remaining case, $|I| < t < l$, write $t' = t - |I|$, $k' = |A'|$ and $l' = |B'|$. It is a simple matter to check that $l' < l$, that $1 \leq t' \leq l'$ and that $k' + l' - D - t' \leq n$. Hence we may apply the induction hypothesis and the evident inequality $\min(\alpha + \beta,\gamma + \delta) \geq \min(\alpha,\gamma) + \min(\beta,\delta)$ to get
\begin{eqnarray*}
\sum_x \min(t,r(A,B,x)) & \geq & \sum_x \min(|I|,r(U,I,x)) + \sum_x \min(t',r(A',B',x)) \\ & = & |U||I| + \sum_x \min(t',r(A',B',x)) \\ & \geq & |U||I| + t'(k' + l' - D - t')\\& = & t(k + l - D - t) + |I|D.\end{eqnarray*}
This is visibly larger than $t\min(n, k + l  - D - t)$, so the induction goes through and the theorem is proved. \endproof\\[11pt]
We will only ever use the following straightforward corollary of Proposition \ref{KP}. If $A,B \subseteq G$ and $K$ is a positive real number write $S_K(A,B) = \{x  :  r(A,B,x) \geq K\}$ for the set of $K$-popular sums in $A + B$.
\begin{corollary}\label{popularsums}
Write $D = D(G)$. Suppose that $A$ and $B$ are subsets of $G$ with cardinalities $k$ and $l$ respectively, and let $K > 0$. Suppose that $\min(k,l) \geq \sqrt{Kn}$. Then
\[ |S_K(A,B)| \; \geq \; \min(n,k + l - D) - 3 \sqrt{Kn}.\]
\end{corollary}
\proof If $K \leq 1$ the result is immediate from Kneser's theorem. Otherwise, for $t \leq \min(k,l)$ we have the inequality
\[ t\min(n,k + l - D - t) \; \leq \; \sum_x \min(t,r(A,B,x)) \; \leq \; Kn + |S_K(A,B)|t.\]
Taking $t = \lfloor \sqrt{Kn}\rfloor$ gives the result (the appearance of the $3$ is because we have taken integer parts).\endproof
\section{Type III groups.} \label{sec7} 
\ni We have determined $\mu(G)$ for all groups except those of type III, whose only prime factors are of the form $6k + 1$. We have also proved Theorem \ref{mainthm2} for such groups. In this section we lay the foundation for a proof of Theorems \ref{mainthm1} and \ref{mainthm2} for type III groups. \\[11pt]
The reader who is only interested in $\mu(G)$ may set $\delta = 0$ throughout the section. This results in some simplification.\\[11pt]
Throughout this section $G$ will be type III, although the same methods can be used to deal with type I and II groups as well (and the calculations are fairly easy in most cases). We begin by setting up some notation. Recall that $\nu(G) = (m-1)/3m$, where $m$ is the exponent of $G$. We are trying to prove that $\mu(G) = \nu(G)$. Let $A \subseteq G$ be a set. Fix a character $\gamma$ for which $\mbox{Re}(\widehat{A}(\gamma))$ is minimal. We call $\gamma$ the \textit{special direction} of $A$. Suppose that $\mbox{ord}(\gamma) = q$, and let $H = \ker \gamma$. Write $H_j = \gamma^{-1}(e^{2\pi ij/q})$ for the cosets of $H$ (so that $H_0 = H$). The indices are to be considered as residues modulo $q$, reflecting the isomorphism $G/H \cong \mathbb{Z}/q\mathbb{Z}$. Write $A_j = A \cap H_j$ and $\alpha_j = |A_j|/|H|$ for the density of $A$ on $H_j$. 
Set $k = (q - 1)/6$, and define the \textit{middle} of $\mathbb{Z}/q\mathbb{Z}$ to be $\{k+1,\dots,5k\}$.\\[11pt]
The next lemma shows how the concept of special direction arises when one is dealing with sets which are nearly sum-free.
\begin{lemma}\label{lem16}
Let $A \subseteq G$ have size $\alpha n$, $\alpha \geq 1/4$, and suppose that there are at most $\delta n^2$ Schur triples in $A$. Then we have
\[ \min_{\gamma} \mbox{\emph{Re}} (\widehat{A}(\gamma)) \; \leq \;  \left(6\delta - \frac{\alpha^2}{1 - \alpha}\right)n.  \]
\end{lemma}
\proof The number of Schur triples in $A$ is exactly $n^{-1}\sum_\gamma |\widehat{A}(\gamma)|^2\widehat{A}(\gamma)$. The contribution from the trivial character $\gamma = 1$ is $\alpha^3 n^2$, and so under the hypotheses of the lemma we have
\[ \sum_{\gamma \neq 1} |\widehat{A}(\gamma)|^2\mbox{Re}(\widehat{A}(\gamma)) \; \leq \; (\delta - \alpha^3)n^3.\]
However by Parseval's identity we have that 
\[ \sum_{\gamma \neq 1} |\widehat{A}(\gamma)|^2 \; = \; (\alpha - \alpha^2)n^2.\] The result follows after a short calculation. \endproof\\[11pt]
The next proposition is the main result of \S \ref{sec7}. For the remainder of the section set $\eta = 2^{-23}$.
\begin{proposition}\label{thm16}
Let $\delta \geq 0$, and suppose that $A \subseteq G$ has at most $\delta n^2$ Schur triples. Let $q = \mbox{\emph{ord}}(\gamma)$ be the order of the special direction in $\Gamma$, let $\kappa = 32\delta^{1/3}q^{2/3}$ and suppose that $\kappa \leq \eta/8q$.  Then either $|A| \leq (\nu(G) - \eta/8)n$ or else $A$ is \emph{essentially middled}, meaning that $\alpha_i \leq 2\kappa$ for all $i$ not in the middle of $\mathbb{Z}/q\mathbb{Z}$.
\end{proposition}
\proof 
The fact that $A$ is nearly sum-free gives a number of inequalities that must be satisfied by the $\alpha_i$. 
These are stated and proved in Lemma \ref{prop9pt1} below, and will be key to our work both in this section and the next.\\[11pt]
Define $\lambda(G)$ to be the \textit{density} of the largest proper subgroup of $G$. Thus $D(G) = \lambda(G)|G|$, where $D(G)$ is the defect as defined in \S \ref{sec6}. Observe also that if $H \leq G$ then $\lambda(H) \leq \lambda(G)$.\\[11pt]
Write $\nu_q = \frac{1}{3}(1 - \frac{1}{q}) = 2k/q$. Since $q$ is the order of the special direction of $A$ it is clear that the exponent $m$ of $G$ satisfies $m \geq q$, and so $\nu(G) \geq \nu_q$.
\begin{lemma}\label{prop9pt1} Let $\delta,A,q$ and $\kappa$ be as in Proposition \ref{thm16}. \begin{enumerate}
\item If $\alpha_l > \kappa$ then for any $j$ we have $\alpha_j + \alpha_{j+l} \;\leq \;1 + \kappa;$
\item For any $i$ we have $\alpha_i + \alpha_{2i} \;\leq \;1 + \kappa;$
\item If $\alpha_u,\alpha_v,\alpha_w > \kappa$, where $u + v = w$, then $\alpha_u + \alpha_v + \alpha_w \;\leq \;1 + \lambda(G) + \kappa;$
\item Suppose that $|A| \geq (\nu(G) - \eta/8)n$. Then \[ q^{-1}\sum_{j=0}^{q-1} \alpha_j \cos \frac{2\pi j}{q} + \frac{\nu_q^2}{1 - \nu_q}\; < \; 3\eta/4.\]
\item Suppose again that $|A| \geq (\nu(G) - \eta/8)n$. Then \[ \sum_{j=0}^{q-1} \alpha_j \; \geq \; (\nu_q - \eta/8) q .\]
\end{enumerate}
\end{lemma}
\proof Observe that the condition $\kappa \leq \eta/8q$, which is one of the assumptions of Proposition \ref{thm16}, implies that $\kappa > \delta^{1/2}q$.\\[11pt] (i) If $\alpha_l > \kappa$ then certainly $\alpha_l > \delta^{1/2}q$, which means that $|A_l| \geq \delta^{1/2}q|H|$. Suppose the result is false, so that 
$\alpha_j + \alpha_{j+l} > 1 + \kappa > 1 + \delta^{1/2}q.$ Thus, for any $x \in A_l$, we have 
\[ \left|(x + A_j) \cap A_{j + l}\right| \; > \; \delta^{1/2}q|H|.\]Hence the number of triples $(x,y,z) \in A_l \times A_j \times A_{l + j}$ with $x + y = z$ is at least $\delta^{1/2}q|H||A_l|$, and so the total number of these triples ober all $l,j$ is greater than $\delta q^2 |H|^2 = \delta n^2$. This is contrary to our assumption.\\[11pt]
(ii) This is immediate from (i).\\[11pt]
(iii) If $\alpha_u,\alpha_v,\alpha_w$ are all greater than $\kappa$ then they are certainly all greater than $5\delta^{1/3}q^{2/3}$. Suppose the result is false, so that $\alpha_u + \alpha_v + \alpha_w > 1 + \lambda(G) + \kappa \geq 1 + \lambda(G) + 5\delta^{1/3}q^{2/3}$.
Observe that $D(H) = \lambda(H)|H| \leq \lambda(G)|H|$. We may apply Corollary \ref{popularsums} with $K =  \delta^{2/3}q^{4/3}|H|$ (clearly $|A_u|,|A_v| \geq \sqrt{K|H|}$), and this gives that
\[
\left|S_K(A_u,A_v)\right| \; \geq \; \min\left(|H|,|A_u| + |A_v| - D(H)\right) - 3\sqrt{K|H|}.\]
We split into the two cases $|H| \geq |A_u| + |A_v| - D(H)$ and $|H| < |A_u| + |A_v| - D(H)$. In the latter case we have
\begin{eqnarray*}
\left|S_K(A_u,A_v) \cap A_w\right| & \geq & |A_w| + |S_K(A_u,A_v)| - |H|\\ & \geq &
5\delta^{1/3}q^{2/3}|H| - 3 \sqrt{K|H|} \\ & \geq & 2\delta^{1/3}q^{2/3}|H|,
\end{eqnarray*}
whereas in the former case we have
\begin{eqnarray*}
\left|S_K(A_u,A_v) \cap A_w \right| & \geq & |A_w| + |S_K(A_u,A_v)| - |H| \\ & \geq & |A_u| + |A_v| + |A_w| - D(H) - |H| - 3\sqrt{K|H|} \\ & \geq & 2\delta^{1/3}q^{2/3}|H|\end{eqnarray*} once again.
Now for any $z \in A_{w} \cap S_K(A_u,A_v)$ there are at least $K$ pairs $(x,y) \in A_u \times A_v$ such that $x + y = z$, and so the total number of triples $(x,y,z) \in A_u \times A_v \times A_{w}$ with $x + y = z$ is at least
\[ 2\delta^{1/3}q^{2/3}|H|K \; > \; \delta q^2|H|^2 \; = \; \delta n^2,\]
a contradiction.\\[11pt]
Recall that for parts (iv) and (v) we are working under the assumption that $|A| \geq (\nu(G) - \eta/8)n$.\\[11pt]
(iv) The fact that $A$ is nearly sum-free implies that $\mbox{Re} (\widehat{A}(\gamma))$ is rather small by Lemma \ref{lem16}. Indeed we have
\[ \mbox{Re}(\widehat{A}(\gamma)) \; = \; |H|\sum_j \alpha_j \cos\frac{2\pi j}{q}\] and so, writing $\alpha = |A|/n$, we have by Lemma \ref{lem16} that
\[ q^{-1}\sum_j \alpha_j \cos\frac{2\pi j}{q} \; \leq \; 6\delta - \frac{\alpha^2}{1 - \alpha}.\] 
Thus
\begin{eqnarray*}
q^{-1}\sum_j \alpha_j \cos \frac{2\pi j}{q} + \frac{\nu(G)^2}{1 - \nu(G)} & \leq & 6\delta  - \frac{\alpha^2}{1 - \alpha} + \frac{\nu(G)^2}{1 - \nu(G)}\\
& \leq & 6\delta + \frac{\eta}{2} \\
& < & \frac{3\eta}{4},
\end{eqnarray*}
the penultimate step following from the fact that $\alpha \geq \nu(G) - \frac{\eta}{8}$. It remains to observe that 
\[ \frac{\nu_q^2}{1 - \nu_q} \; \leq \; \frac{\nu(G)^2}{1 - \nu(G)},\] a consequence of the inequality $\nu(G) \geq \nu_q$.\\[11pt]
(v) Immediate from the fact that $\nu(G) \geq \nu_q$.\endproof\\[11pt]
We will need to do several calculations with the inequalities of Lemma \ref{prop9pt1}, and for that reason it will be convenient to have them in an easy-to-use form. 
\begin{definition}\label{defbeta} Define numbers $\beta_i \in [0,1]$, $i = 0,\dots,q-1$, by
$\beta_i  =  0$ if $\alpha_i \leq \kappa$ and \[ \beta_i \; = \; \frac{\alpha_i - \kappa}{1 + \kappa} \] otherwise.
\end{definition}
\ni An important property of the $\beta_i$ is their relation to the $\alpha_i$, which is that these two sets of numbers are rather close. Indeed it is easy to confirm that 
\begin{equation}\label{eq1979} |\beta_i - \alpha_i| \; \leq \; 2\kappa. \end{equation}
The next lemma details the inequalities satisfied by the $\beta_i$. 
\begin{lemma}\label{betaineqs}Let $\delta \geq 0$, and suppose that $A$ has at most $\delta n^2$ Schur triples. \begin{enumerate}
\item If $\beta_l > 0$ then for any $j$ we have $\beta_j + \beta_{j+l} \;\leq \;1;$
\item For any $i$ we have $\beta_i + \beta_{2i} \;\leq \;1;$
\item If $\beta_u,\beta_v,\beta_w > 0$, where $u + v = w$, then $\beta_u + \beta_v + \beta_w \;\leq \;1 + \lambda(G);$
\item If $\beta_{2t} > 0$ then $2\beta_t + \beta_{2t} \;\leq \;1 + \lambda(G)$.
\item Suppose that $|A| \geq (\nu(G) - \eta/8)n$. Then \[ q^{-1}\sum_{j=0}^{q-1} \beta_j \cos \frac{2\pi j}{q} + \frac{\nu_q^2}{1 - \nu_q}\; < \; \eta.\]
\item Suppose again that $|A| \geq (\nu(G) - \eta/8)n$. Then \[ \sum_{j=0}^{q-1} \beta_j \; \geq \; (\nu_q - \eta) q \; = \; 2k - \eta q.\]
\end{enumerate}
\end{lemma}
\proof (i),(ii) and (iii) follow quickly from Lemma \ref{prop9pt1}. (iv) follows from (iii) if $\beta_t > 0$, and is immediate if $\beta_t = 0$. For (v) and (vi), we use \eqref{eq1979} and recall that $\kappa \leq \eta/8q$. \endproof\\[11pt]
We now return in earnest to the proof of Proposition \ref{thm16}. Our strategy will be this: assume that $A$ is not essentially middled and at the same time that $|A| > (\nu(G) - \eta/8)n$. We will obtain a contradiction using the inequalities of Lemma \ref{betaineqs}. Under these assumptions, (vi) of Lemma \ref{betaineqs} certainly holds. Also (by reflecting $A$ if necessary) we may assume that $\beta_l > 0$ for some $l \in \{0,1,\dots,k\}$, which we assume to be minimal. Lemma \ref{betaineqs} (i) then tells us that $\beta_j + \beta_{j + l} \leq 1$ for all $j$.\\[11pt] Write $M$ for the quantity in Lemma \ref{betaineqs} (v), that is 
\[ M \; = \; q^{-1}\sum_{j=0}^{q-1} \beta_j \cos\frac{2\pi j}{q} + \frac{\nu_q^2}{1 - \nu_q}.\] Our aim will be to contradict Lemma \ref{betaineqs} (v). It is important to note that the preceding assumptions and notation will be used for the remainder of this section.
\begin{lemma}\label{prop2} Let $\delta, A, q$ and $\kappa$ be as in Proposition \ref{thm16}. Assume that $|A| > (\nu(G) - \eta/8)n$, that $A$ is not essentially middled, and suppose that that $l \in \{ 0,1,\dots,k\}$ satisfies $0 = \beta_0 = \dots = \beta_{l-1} < \beta_l$. Then we have 
\begin{equation}\label{ii1} M \; \geq \; -\frac{\sin (4 k\pi/q)}{2q \sin (\pi/q)\cos (\pi l/q)} + \frac{\nu_q^2}{1  - \nu_q}- 6\eta\end{equation} if $l$ is even and
\begin{equation}\label{ii2} M \; \geq \; -\frac{\sin((4k-1) \pi /q)}{2q\sin(\pi/q)\cos(\pi l/q)} -  \frac{\cos (4k\pi/q)}{2q \cos(\pi l /q)} + \frac{\nu_q^2}{1  - \nu_q} - 6\eta\end{equation} if $l$ is odd. 
\end{lemma}
\proof Naturally, we use the inequalities $\beta_j + \beta_{j + l} \leq 1$. We have
\begin{equation}\label{eq77}
\cos \frac{\pi l}{q}\sum_j \beta_j \cos\frac{2\pi j}{q} \; = \;  \frac{1}{2}\sum_j (\beta_j + \beta_{j+l}) \cos\frac{2j+l}{q}\pi \; = \; \sum_j \gamma_j \cos\frac{2j+l}{q}\pi, \end{equation}
where $\gamma_j := (\beta_j + \beta_{j+l})/2$.
Write $E(s)$ for the minimum value of $\sum_j \gamma_j \cos \frac{2j + l}{q}\pi$ subject to the constraints $0 \leq \gamma_j \leq 1/2$ and $\sum_j \gamma_j \geq 2k - sq$. By virtue of Lemma \ref{betaineqs} (vi) we are, of course, interested in $E(\eta)$. The following lemma allows us to concentrate on the somewhat more appealing quantity $E(0)$.
\begin{lemma}\label{tedious-lemma}$E(\eta) \geq E(0) - 3\eta q$.
\end{lemma}
\proof Write $c_j = \cos \frac{2j + l}{q}\pi$. Suppose that the numbers $\gamma_j$ satisfy the constraints $0 \leq \gamma_j \leq 1/2$ and $\sum_j \gamma_j \geq 2k - \eta q$, and that $\sum_j \gamma_j c_j = E(\eta)$. Define
\[ \gamma'_j := \lambda \gamma_j + \frac{1}{2}(1 - \lambda),\]
where 
\[ 1/\lambda := 1 + \frac{2\eta q}{q - 4k}.\]
Since $\lambda < 1$, the numbers $\gamma'_j$ satisfy $0 \leq \gamma'_j \leq 1/2$. It is easy to see that $\sum_j \gamma'_j \geq 2k$, and so we have
\[ E(0) \leq \sum_j \gamma'_j c_j \leq E(\eta) + q\max_j |\gamma_j - \gamma'_j|.\]
Writing $x = 2\eta q/(q - 4k)$, and observing that $x \leq 6\eta$, the proof is concluded by the noting the chain of inequalities
\begin{equation} \label{boxeq}|\gamma'_j - \gamma_j| \leq \frac{1}{2}|1 - \lambda| = \frac{x}{2(1 + x)} \leq x/2 \leq 3\eta.\end{equation}
Returning now to the proof of Lemma \ref{prop2}, write $u := (q-1)/2$ and suppose first that $l = 2m$ is even. Here the minimum in question is obtained by setting $\gamma_j = 1/2$ for $j + m = u - 2k + 1,\dots, u,u + 1,\dots,u + 2k$ and $\gamma_j = 0$ otherwise.
Observe that $2k \geq (q + 1)/4$, ensuring that in the extremal configuration all the negative $\cos((2j + l)\pi/q)$ have the largest possible weight $\gamma_j = 1/2$. Thus the best strategy is to put as few weights on the positive values as possible. A precise computation is possible and yields
\begin{equation}\label{eq16} E(0) \; = \; -\frac{\sin(4k \pi/q)}{2 \sin(\pi/q)} \qquad\qquad \mbox{($l$ even)}.\end{equation} 
The case $l = 2m + 1$ is very similar; the minimum is also obtained by setting $\gamma_j = 1/2$ when  $j + m = u - 2k + 1,\dots,u,u+1,\dots,u + 2k$ and $\gamma_j = 0$ otherwise. An exact evaluation is again possible and leads to 
\begin{equation}\label{eq17} E(0) \; = \;  -\frac{\sin((4k-1)\pi /q)}{2\sin(\pi/q)} - \frac{1}{2}\cos(4\pi k/q) \qquad \mbox{($l$ odd)}.\end{equation}
The desired inequalities \eqref{ii1} and \eqref{ii2} follow from Lemma \ref{tedious-lemma}, equations \eqref{eq77}, \eqref{eq16} and \eqref{eq17} and the fact that $\cos(\pi l/q) \geq\cos(\pi k/q) > 1/2$.\endproof
\begin{lemma}\label{lem34}  Let $\delta, A, q$ and $\kappa$ be as in Proposition \ref{thm16}. Assume that $|A| > (\nu(G) - \eta/8)n$, that $A$ is not essentially middled, and suppose that that $l \in \{ 0,1,\dots,k\}$ satisfies $0 = \beta_0 = \dots = \beta_{l-1} < \beta_l$.
 Then we have the bound $M \geq \eta$ \emph{(}implying Lemma \ref{thm16}\emph{)} in the following cases:\begin{enumerate}
\item $l \leq k - 2$ except for $q = 7,13,19$;
\item $l = k - 1$ except for $q = 7,13,19,31$;
\item $l = k$ except for $q = 7,13,19,31,37,43,49,61,67,73$.
\end{enumerate}
\end{lemma}
\proof The authors proved this with the aid of a computer. We confine ourselves here to a few remarks which would enable the extremely keen reader to reproduce our computations.\\[11pt] For large $q$, we have the approximations $2k/q \approx 1/3$ and $q\sin(\pi/q) \approx \pi$. Since $l \leq k$, we have $\cos(\pi l/q) \geq \cos(\pi / 6) = \sqrt{3}/2$. Substituting into \eqref{ii1} shows that $M$ is at least $\frac{1}{6} - \frac{1}{2\pi} + \epsilon_q - 6\eta$, where $\epsilon_q \rightarrow 0$ as $q \rightarrow \infty$, uniformly in $l$. This is greater than $\eta$ for $q$ sufficiently large. By estimating $\epsilon_q$ precisely using simple (but rather tedious) calculus one may verify that a suitable notion of sufficiently large is $q > 1000$. Equation \eqref{ii2} may be treated in exactly the same way. The remaining pairs $(l,q)$, where $q < 1000$, can be checked individually on a computer and Lemma \ref{lem34} is what results.\endproof\\[11pt]
There are several cases not covered by Lemma \ref{lem34}. The great majority can be dealt with by using the additional linear relations $\beta_i + \beta_{2i} \leq 1$. The problem of minimising $M$ subject to these constraints, the constraints $\beta_i + \beta_{i + l} \leq 1$, $\beta_i \geq 0$ and the inequality $\sum \beta_i \geq (\nu_q - \eta)q$ is a standard linear programming problem, but we have not found a convenient way of dealing with it for general $q,l$. For any specific values however one may consider the \textit{dual} problem in the sense of linear programming theory. 
Recall that if we have a \textit{primal} problem 
\begin{equation}\label{x1} \mbox{maximise $r^T x$ subject to $Bx \leq c$, $x \geq 0$}\end{equation}
then we may associate to it a dual problem 
\begin{equation}\label{x2} \mbox{minimise $y^T c$ subject to $y^T B \geq r^T$, $y \geq 0$}.\end{equation}
If the solutions to these problems are $S$ and $S'$ respectively then one has $S = S'$. The equality is a reasonably deep theorem, but it is easy to see that $S \leq S'$, and this is all we need later on. Indeed, for any $x$ and $y$ satisfying $Bx \leq c$ and $y^TB \geq r^T$ we have $r^T x \leq y^T B x \leq y^T c$. That is, for any values of $y$ that satisfy $y^T B \geq r^T$, $y \geq 0$ the quantity $y^T c$ provides an upper bound for the primal problem.\\[11pt]
Our problem may be cast in the form \eqref{x1} simply by multiplying the constraint $\sum \beta_i \geq (\nu_q - \eta)q$ by $-1$. In the dual problem, the constraints $\beta_i + \beta_{2i} \leq 1$, $\beta_i + \beta_{i+l} \leq 1$ become variables $\lambda_i,\mu_i$ respectively, and the constraint $-\sum \beta_i \leq -(\nu_q - \eta)q$ becomes a variable $\tau$. The dual problem is then to maximise
\begin{equation} \label{dual}(\nu_q - \eta)q - \lambda_1 - \dots - \lambda_q - \mu_1 - \dots - \mu_q + \frac{\nu_q^2}{1  - \nu_q}\end{equation} subject to the relations $\lambda_i \geq 0$, $\mu_i \geq 0$, $\tau \geq 0$ and  
\begin{equation}\label{constraints} \tau - \lambda_{j} - \lambda_{j - l} - \mu_j - \mu_{j/2} \; \leq \; \cos\frac{2\pi j}{q}\end{equation} for $j = 1,\dots,q$ (where addition and division by 2 are of course taken modulo $q$). Values are provided for $\tau,\lambda_j,\mu_j$ in the lists below, giving lower bounds for $M$.  Of course, these values were found using a simplex algorithm on a computer; however because of the duality result $S \leq S'$ the reader need not concern herself with the numerical accuracy of our routines. Any value which does not appear explicitly is assigned the value $0$.\\[11pt]
\noindent $(q,l) = (73,12)$. $M \geq 0.01$. 
$\scs\lambda_{13} = 0.072$, $\scs\lambda_{15} = 0.237$, $\scs\lambda_{18} = 0.491$, $\scs\lambda_{19} = 0.465$, $\scs\lambda_{20} = 0.663$, $\scs\lambda_{23} = 0.910$, $\scs\lambda_{24} = 0.757$, $\scs\lambda_{25} = 0.458$, $\scs\lambda_{26} = 1.127$, $\scs\lambda_{27} = 0.707$, $\scs\lambda_{28} = 0.723$, $\scs\lambda_{29} = 1.029$, $\scs\lambda_{30} = 0.870$, $\scs\lambda_{31} = 0.936$, $\scs\lambda_{32} = 0.456$, $\scs\lambda_{33} = 0.999$, $\scs\lambda_{34} = 1.085$, $\scs\lambda_{35} = 0.595$, $\scs\lambda_{36} = 0.755$, $\scs\lambda_{37} = 0.988$, $\scs\lambda_{38} = 0.378$, $\scs\lambda_{39} = 0.783$, $\scs\lambda_{40} = 0.745$, $\scs\lambda_{41} = 0.194$, $\scs\lambda_{42} = 0.324$, $\scs\lambda_{43} = 0.425$, $\scs\lambda_{44} = 0.026$, $\scs\lambda_{45} = 0.104$, $\scs\lambda_{48} = 0.076$, $\scs\lambda_{57} = 0.217$, $\scs\mu_{13} = 0.005$, $\scs\mu_{14} = 0.155$, $\scs\mu_{16} = 0.320$, $\scs\mu_{17} = 0.405$, $\scs\mu_{21} = 0.210$, $\scs\mu_{22} = 0.830$, $\scs\mu_{24} = 0.231$, $\scs\mu_{25} = 0.533$, $\scs\mu_{27} = 0.254$, $\scs\mu_{28} = 0.380$, $\scs\mu_{29} = 0.237$, $\scs\mu_{46} = 0.112$, $\scs\mu_{47} = 0.537$, $\scs\mu_{51} = 0.047$, $\scs\mu_{52} = 0.002$, $\scs\mu_{53} = 0.469$, $\scs\mu_{55} = 0.067,$
$\scs\mu_{59} = 0.155$, $\scs\tau = 0.512$.\\[11pt]
$(q,l) = (67,11)$. $M \geq 0.01$.
$\scs\lambda_{13} = 0.087$, $\scs\lambda_{14} = 0.176$, $\scs\lambda_{15} = 0.268$, $\scs\lambda_{16} = 0.361$, $\scs\lambda_{17} = 0.455$, $\scs\lambda_{18} = 0.364$, $\scs\lambda_{19} = 0.452$, $\scs\lambda_{20} = 0.632$, $\scs\lambda_{21} = 0.637$, $\scs\lambda_{22} = 0.675$, $\scs\lambda_{23} = 0.756$, $\scs\lambda_{24} = 0.974$, $\scs\lambda_{25} = 0.583$, $\scs\lambda_{26} = 0.927$, $\scs\lambda_{27} = 0.796$, $\scs\lambda_{28} = 0.847$, $\scs\lambda_{29} = 0.792$, $\scs\lambda_{30} = 0.927$, $\scs\lambda_{31} = 0.772$, $\scs\lambda_{32} = 0.785$, $\scs\lambda_{33} = 0.756$, $\scs\lambda_{34} = 0.675$, $\scs\lambda_{35} = 0.448$, $\scs\lambda_{36} = 0.637$, $\scs\lambda_{37} = 0.452$, $\scs\lambda_{38} = 0.548$, $\scs\lambda_{39} = 0.455$, $\scs\lambda_{40} = 0.361$, $\scs\lambda_{41} = 0.268$, $\scs\lambda_{42} = 0.176$, $\scs\lambda_{43} = 0.087$, $\scs\mu_{18} = 0.185$,
$\scs\mu_{20} = 0.099$, $\scs\mu_{21} = 0.183$, $\scs\mu_{23} = 0.229$, $\scs\mu_{43} = 0.189$, $\scs\mu_{45} = 0.229$,
$\scs\mu_{46} = 0.372$, $\scs\mu_{47} = 0.095$, $\scs\mu_{48} = 0.189$, $\scs\tau = 0.430$.\\[11pt]
$(q,l) = (61,10)$. $M \geq 0.01$.
$\scs\lambda_{12} = 0.096$, $\scs\lambda_{13} = 0.195$, $\scs\lambda_{14} = 0.296$, $\scs\lambda_{15} = 0.399$, $\scs\lambda_{16} = 0.502$, $\scs\lambda_{17} = 0.448$, $\scs\lambda_{18} = 0.647$, $\scs\lambda_{19} = 0.599$, $\scs\lambda_{20} = 0.668$, $\scs\lambda_{21} = 0.756$, $\scs\lambda_{22} = 0.814$, $\scs\lambda_{23} = 0.740$,
$\scs\lambda_{24} = 0.913$, $\scs\lambda_{25} = 0.818$, $\scs\lambda_{26} = 0.818$, $\scs\lambda_{27} = 0.913$,
$\scs\lambda_{28} = 0.745$, $\scs\lambda_{29} = 0.814$, $\scs\lambda_{30} = 0.756$, $\scs\lambda_{31} = 0.668$,
$\scs\lambda_{32} = 0.599$, $\scs\lambda_{33} = 0.652$, $\scs\lambda_{34} = 0.448$, $\scs\lambda_{35} = 0.502$,
$\scs\lambda_{36} = 0.394$, $\scs\lambda_{37} = 0.296$, $\scs\lambda_{38} = 0.195$, $\scs\lambda_{39} = 0.096$,
$\scs\mu_{18} = 0.057$, $\scs\mu_{19} = 0.202$, $\scs\mu_{20} = 0.227$, $\scs\mu_{22} = 0.157$, $\scs\mu_{23} = 0.005$,
$\scs\mu_{39} = 0.157$, $\scs\mu_{41} = 0.227$, $\scs\mu_{42} = 0.202$, $\scs\mu_{43} = 0.053$, $\scs\tau = 0.423$.\\[11pt]
$(q,l) = (49,8)$. $M \geq 0.01$.
$\scs\lambda_{10} = 0.121$, $\scs\lambda_{11} = 0.246$, $\scs\lambda_{12} = 0.373$, $\scs\lambda_{13} = 0.278$, $\scs\lambda_{14} = 0.628$, $\scs\lambda_{15} = 0.430$, $\scs\lambda_{16} = 0.710$, $\scs\lambda_{17} = 0.536$, $\scs\lambda_{18} = 0.892$, $\scs\lambda_{19} = 0.885$, $\scs\lambda_{20} = 0.871$, $\scs\lambda_{21} = 0.871$, $\scs\lambda_{22} = 0.727$, $\scs\lambda_{23} = 0.957$, $\scs\lambda_{24} = 0.694$, $\scs\lambda_{25} = 0.868$, $\scs\lambda_{26} = 0.272$, $\scs\lambda_{27} = 0.470$, $\scs\lambda_{28} = 0.436$, $\scs\lambda_{29} = 0.373$, $\scs\lambda_{30} = 0.246$, $\scs\lambda_{31} = 0.121$, $\scs\mu_{13} = 0.223$, $\scs\mu_{15} = 0.195$, $\scs\mu_{16} = 0.158$, $\scs\mu_{17} = 0.442$, $\scs\mu_{18} = 0.066$, $\scs\mu_{32} = 0.127$, $\scs\mu_{34} = 0.037$, $\scs\mu_{35} = 0.158$, $\scs\tau = 0.404$.\\[11pt]
$(q,l) = (43,7)$. $M \geq 0.01$.
$\scs\lambda_{9} = 0.139$, $\scs\lambda_{10} = 0.282$, $\scs\lambda_{11} = 0.428$, $\scs\lambda_{12} = 0.573$, $\scs\lambda_{13} = 0.421$, $\scs\lambda_{14} = 0.632$, $\scs\lambda_{15} = 0.612$, $\scs\lambda_{16} = 0.948$,
$\scs\lambda_{17} = 0.898$, $\scs\lambda_{18} = 0.836$, $\scs\lambda_{19} = 0.608$, $\scs\lambda_{20} = 0.948$,
$\scs\lambda_{21} = 0.757$, $\scs\lambda_{22} = 0.778$, $\scs\lambda_{23} = 0.421$, $\scs\lambda_{24} = 0.428$,
$\scs\lambda_{25} = 0.428$, $\scs\lambda_{26} = 0.282$, $\scs\lambda_{27} = 0.139$,
$\scs\mu_{13} = 0.294$, $\scs\mu_{14} = 0.217$, $\scs\mu_{15} = 0.291$, $\scs\mu_{29} = 0.072$, $\scs\mu_{30} = 0.004$,
$\scs\mu_{31} = 0.146$, $\scs\tau = 0.391$.\\[11pt]
$(q,l) = (37,6)$. $M \geq 0.01$.
$\scs\lambda_{8} = 0.163$, $\scs\lambda_{9} = 0.331$, $\scs\lambda_{10} = 0.414$, $\scs\lambda_{11} = 0.408$,
$\scs\lambda_{12} = 0.742$, $\scs\lambda_{13} = 0.628$, $\scs\lambda_{14} = 0.933$, $\scs\lambda_{15} = 0.785$,
$\scs\lambda_{16} = 0.871$, $\scs\lambda_{17} = 0.933$, $\scs\lambda_{18} = 0.628$, $\scs\lambda_{19} = 0.742$,
$\scs\lambda_{20} = 0.322$, $\scs\lambda_{21} = 0.500$, $\scs\lambda_{22} = 0.331$, $\scs\lambda_{23} = 0.163$,
$\scs\mu_{10} = 0.087$, $\scs\mu_{12} = 0.082$, $\scs\mu_{13} = 0.258$, $\scs\mu_{24} = 0.258$, $\scs\mu_{25} = 0.082$,
$\scs\mu_{26} = 0.087$, $\scs\tau = 0.373$.\\[11pt]
$(q,l) =(31,5)$. $M \geq 0.005$.
$\scs\lambda_{7} = 0.196$, $\scs\lambda_{8} = 0.398$, $\scs\lambda_{9} = 0.392$, $\scs\lambda_{10} = 0.689$,
$\scs\lambda_{11} = 0.654$, $\scs\lambda_{12} = 0.911$, $\scs\lambda_{13} = 0.824$, $\scs\lambda_{14} = 0.911$,
$\scs\lambda_{15} = 0.654$, $\scs\lambda_{16} = 0.689$, $\scs\lambda_{17} = 0.392$, $\scs\lambda_{18} = 0.398$,
$\scs\lambda_{19} = 0.196$, $\scs\mu_{10} = 0.100$, $\scs\mu_{11} = 0.207$, $\scs\mu_{20} = 0.207$, $\scs\mu_{21} = 0.100$, $\scs\tau = 0.347$.\\[11pt]
$(q,l) =(31,4)$. $M \geq 0.005$.
$\scs\lambda_{8} = 0.081$, $\scs\lambda_{10} = 0.254$, $\scs\lambda_{11} = 0.635$, $\scs\lambda_{12} = 1.026$,
$\scs\lambda_{13} = 0.732$, $\scs\lambda_{14} = 0.853$, $\scs\lambda_{15} = 0.709$, $\scs\lambda_{17} = 0.571$,
$\scs\lambda_{19} = 0.398$, $\scs\lambda_{20} = 0.196$, $\scs\mu_{7} = 0.196$, $\scs\mu_{8} = 0.318$, $\scs\mu_{9} = 0.370$,
$\scs\mu_{10} = 0.535$, $\scs\mu_{11} = 0.108$, $\scs\mu_{20} = 0.229$, $\scs\mu_{21} = 0.218$, $\scs\mu_{22} = 0.491$,
$\scs\tau = 0.347$.\\[11pt]
For $q = 7,13,19$ it is necessary to consider inequalities (iii) and (iv) of Proposition \ref{betaineqs}. Since $G$ is type III we have $\lambda(G) \leq 1/7$. Now, however, the problem is not strictly linear because the new constraints only hold if certain of the variables involved are known to be strictly positive. For this reason it is necessary to split into cases. In what follows we will be concerned with the minimisation of $M$ subject to the constraints $\beta_j + \beta_{j+l} \leq 1$, $\beta_j + \beta_{2j} \leq 1$, $\sum \beta_i \geq (\nu_q - \eta)q$ and certain other linear constraints $C_1,C_2,\dots$ which will be listed in each case. As before we will consider the dual problem. This problem will have variables $\lambda_1,\dots,\lambda_q$, $\mu_1,\dots,\mu_q$, $\tau$ and $\theta_1,\theta_2,\dots$, the constraint $C_j$ giving rise to the dual variable $\theta_j$.\\[11pt]
It saves some time to observe that there is no need to treat the case $l = 0$ separately. Indeed the inequalities $\beta_j + \beta_{j+l} \leq 1$ for $l = 0$ imply those for all other $l$, and (as the reader may care to check) we do not make any further use of the fact that $\beta_l > 0$.\\[11pt]
\textbf{The prime 19.} $\nu_{19} = 6/19$. To deal with the case $q = 19$ we only need, for any $l$, the single extra constraint $2\beta_7 + \beta_{14} \leq 4/3$. However we only know this in the case $\beta_{14} > 0$ and so we must deal separately with the case $\beta_{14} = 0$ for each $l$. Thus we must solve the following six linear problems.\\[11pt]
\noindent $(q,l) = (19,3)$, $C_1 : \beta_{14} \leq 0$. $M \geq 0.01$.
$\scs\lambda_4 = 0.302$, $\scs\lambda_5 = 0.302$, $\scs\lambda_6 = 0.629$, $\scs\lambda_7 = 0.590$, $\scs\lambda_8 1.125$, $\scs\lambda_9 0.905$, $\scs\lambda_{10} = 0.616$, $\scs\mu_5 = 0.329$, $\scs\mu_6 0.320$, $\scs\mu_{13} = 0.334,
\scs\mu_{15} = 0.302$, $\scs\theta_1 = 0.630$, $\scs\tau = 0.546$.\\[11pt]
\noindent $(q,l) = (19,3)$, $C_1 : 2\beta_7 + \beta_{14} \leq 4/3$. $M \geq 0.004$.
$ \scs\lambda_{6} = 0.6600$, $\scs\lambda_{8} = 1.1250$, $\scs\lambda_{9} = 0.8734$,
$\scs\lambda_{10} = 0.9487$, $\scs\lambda_{11} = 0.0175$, $\scs\lambda_{12} = 0.0175$, $\scs\mu_{4} = 0.3015$,
$\scs\mu_{5} = 0.5847$, $\scs\mu_{6} = 0.2887$, $\scs\mu_{12} = 0.0449$, $\scs\mu_{15} = 0.2841$, $\scs\theta_1 = 0.6122$,
$\scs\tau =  0.5469$.\\[11pt]
\noindent $(q,l) = (19,2)$, $C_1 : \beta_{14} \leq 0$. $M \geq 0.01$.
$\scs\lambda_4 = 0.302$, $\scs\lambda_6 = 0.648$, $\scs\lambda_7 = 1.225$, $\scs\lambda_8 = 0.780$,
$\scs\lambda_9 = 0.310$, $\scs\lambda_{10} = 0.675$, $\scs\lambda_{11} = 0.883$, $\scs\lambda_{13} = 0.067$,
$\scs\mu_5 = 0.080$, $\scs\mu_{12} = 0.550$, $\scs\mu_{15} = 0.236$, $\scs\theta_1 = 0.630$, $\scs\tau = 0.547$.\\[11pt]
\noindent $(q,l) = (19,2)$, $C_1 : 2\beta_7 + \beta_{14} \leq 4/3$. $M \geq 0.002$.
$\scs\lambda_{4} = 0.3015$, $\scs\lambda_{6} = 0.6472$, $\scs\lambda_{7} = 0.4912$, $\scs\lambda_{8} = 0.7793$, $\scs\lambda_{9} = 0.7793$, $\scs\lambda_{10} = 0.6744$, $\scs\lambda_{11} = 0.6472$, $\scs\lambda_{13} = 0.3015$, $\scs\mu_{5} = 0.0797$, $\scs\mu_{12} = 0.5499$, $\scs\mu_{14} =  0.2630$, $\scs\theta_1 = 0.3666$,
$\scs\tau = 0.5469$.\\[11pt]
\noindent $(q,l) = (19,1)$, $C_1 : \beta_{14} \leq 0$. $M \geq 0.01$.
$\scs\lambda_4 = 0.302$, $\scs\lambda_5 = 0.329$, $\scs\lambda_6 = 0.621$, $\scs\lambda_7 = 0.591$,
$\scs\lambda_8 = 0.837$, $\scs\lambda_9 = 0.698$, $\scs\lambda_{10} = 0.837$, $\scs\lambda_{11} = 0.289$,
$\scs\lambda_{12} = 0.936$, $\scs\mu_{13} = 0.014$, $\scs\mu_{15} = 0.302$, $\scs\theta_1 = 0.630$, $\scs\tau = 0.546$.\\[11pt]
\noindent $(q,l) = (19,1)$, $C_1 : 2\beta_7 + \beta_{14} \leq 4/3$. $M \geq 0.004$.
$\scs\lambda_{5} = 0.6296$, $\scs\lambda_{6} = 0.0175$, $\scs\lambda_{8} = 1.1250$, $\scs\lambda_{9} = 0.4084$, $\scs\lambda_{10} = 1.1250$, $\scs\lambda_{12} = 0.9226$, $\scs\lambda_{13} = 0.0262$, $\scs\mu_{4} = 0.3015$,
$\scs\mu_{6} = 0.3018$, $\scs\mu_{15} = 0.3015$, $\scs\theta_1 = 0.6035$, $\scs\tau = 0.5469$.\\[11pt]
\textbf{The prime 13.} $\nu_{13} = 4/13$. Proposition \ref{betaineqs} tells us that either $\beta_3 = 0$, $\beta_{10} = 0$ or both of the relations $2\beta_5 + \beta_{10} \leq 8/7$ and $2\beta_8 + \beta_3 \leq 8/7$ are true.  We will show that any of these conditions suffices by itself to show $M$ positive. Furthermore by symmetry we only need deal with one of the cases $\beta_3 = 0$ and $\beta_{10} = 0$. Thus for each $l$ there are two linear problems to consider, giving the following four problems in all.\\[11pt]
$(q,l) = (13,2)$, $C_1 : \beta_{10} \leq 0$. $M \geq 0.003$.
$\scs\lambda_3 = 0.4476$, $\scs\lambda_4 = 0.5726$, $\scs\lambda_5 = 0.7427$, $\scs\lambda_6 = 0.9665$, $\scs\lambda_7 = 0.7964$, $\scs\mu_4 = 0.3502$, $\scs\mu_9 = 0.1264$, $\scs\theta_1 = 0.4476$, $\scs\tau = 0.5680$.\\[11pt]
$(q,l) = (13,2)$, $C_1 : 2\beta_5 + \beta_{10} \leq 8/7$, $C_2: 2\beta_8 + \beta_3 \leq 8/7$. $M \geq 0.007$.
$\scs\lambda_{4} = 0.9227$, $\scs\lambda_{5} = 0.5514$, $\scs\lambda_{6} = 0.5514$, $\scs\lambda_{7} = 0.9227$, $\scs\lambda_{10} = 0.0650$, $\scs\mu_3 = 0.0650$, $\scs\theta_1 = 0.3826$, $\scs\theta_2 = 0.3826$,
$\scs\tau = 0.5680$.\\[11pt]
$(q,l) = (13,1)$, $C_1 : \beta_{10} \leq 0$. $M \geq 0.003$.
$\scs\lambda_3 = 0.4476$, $\scs\lambda_4 = 0.4752$, $\scs\lambda_5 = 0.6177$, $\scs\lambda_6 = 0.9214$, $\scs\lambda_7 = 0.6177$, $\scs\lambda_8 = 0.6990$, $\scs\mu_9 = 0.2238$, $\scs\theta_1 = 0.4476$, $\scs\tau = 0.5680$.\\[11pt]
$(q,l) = (13,1)$, $C_1 : 2\beta_5 + \beta_{10} \leq 8/7$, $C_2: 2\beta_8 + \beta_3 \leq 8/7$.
$M \geq 0.001$.
$\scs\lambda_{3} = 0.3341$, $\scs\lambda_{4} = 0.4216$, $\scs\lambda_{6} = 1.5391$, $\scs\lambda_{8} = 0.9227$, $\scs\mu_4 = 0.1671$, $\scs\theta_1 = 0.4476$, $\scs\theta_2 = 0.1135$, $\scs\tau = 0.5680$.\\[11pt]
\noindent\textbf{The prime 7.} $\nu_7 = 2/7$. Finally we come to what is, in some sense, the trickiest case. Our workload is reduced by only having to consider the single case $l = 1$. 
Since $\nu_7 = 2/7$ we have $\sum \beta_i \geq 2 - \eta$. We begin by reducing to the case $\beta_2,\beta_3,\beta_4,\beta_5 > 0$. Indeed consider the linear problem with the 15 usual constraints and the condition $C_1 : \beta_3 \leq 0$. Then we have $M \geq 0.04$ as is shown by the dual values
$\lambda_4 = \lambda_9 = 0.679$, $\theta_1 = 1.357$, $\tau = 0.456$. Similarly if $\beta_2 = 0$ we have $M \geq 0.01$, with the dual values $\lambda_3 = 1.102$, $\lambda_4 = 0.424$, $\mu_5 = 0.424$, $\theta_1 = 0.847$, $\tau = 0.623$. The cases $\beta_4 = 0$ and $\beta_5 = 0$ follow by symmetry.\\[11pt]
Suppose then that $\beta_2,\beta_3,\beta_4,\beta_5 > 0$, but that $\beta_1 = \beta_6 = 0$. Since we are assuming that $A$ is not middled we must have $\beta_0 > 0$. But this would imply
\begin{eqnarray*}
4\left(\beta_0 + \sum_{i=2}^5 \beta_i\right) & = & \sum_{i=2}^5 (\beta_0 + 2\beta_i) + (2\beta_2 + \beta_4) + (2\beta_5 + \beta_3) + (\beta_3 + \beta_4) \\ & \leq & \textstyle\frac{32}{7} + \frac{8}{7} + \frac{8}7 + 1 \; = \; 4\left(2 - \frac{1}{28}\right),
\end{eqnarray*}
a contradiction.
It follows that either $\beta_1 > 0$ or $\beta_6 > 0$. Suppose that $\beta_1 = 0$. Then $\beta_6 > 0$, and we solve another linear problem with the 15 usual constraints together with $C_1 : \beta_1 \leq 0$, $C_2 : \beta_4 + \beta_5 + \beta_6 \leq 8/7$, $C_3 : 2\beta_3 + \beta_6 \leq 8/7$, $C_4: 2\beta_5 + \beta_3 \leq 8/7$, $C_5: 2\beta_2 + \beta_4\leq 8/7$. We have $M \geq 0.002$, as witnessed by dual values
$\lambda_{3} = 1.0746$, $\theta_1 = 0.3766$, $\theta_2 = 0.1614$, $\theta_3 = 0.5037$, $\theta_4 = 0.6113$, $\theta_5 = 0.2152$, $\tau = 1.0000$.\\[11pt]
Finally we deal with the case $\beta_1,\beta_2,\dots,\beta_6 > 0$. Once again we take the 15 usual constraints, plus
$C_1 : 2\beta_3 + \beta_6 \leq 8/7$, $C_2 : \beta_4 + \beta_5 + \beta_2 \leq 8/7$, $C_3 : \beta_2 + \beta_4 + \beta_6 \leq 8/7$, $C_4 : \beta_2 + \beta_3 + \beta_5 \leq 8/7$, $C_5 : \beta_1 + \beta_3 + \beta_4 \leq 8/7$, $C_6 : 2\beta_5 + \beta_3 \leq 8/7$, $C_7: 2\beta_2 + \beta_4 \leq 8/7$. With these conditions one has $M \geq 0.0005$ as shown by the values
$\lambda_3 = 0.53669$, $\theta_1 = 0.37652$, $\theta_4 = 0.98778$, $\theta_6 = 0.23476$, $\theta_7 = 0.37652$, $\tau = 1.0000$.\\[11pt]
The proof of Proposition \ref{thm16} is at long last complete.\endproof
\section{Sum-free sets in type III groups.} 
\ni The next result completes the proof of Theorem \ref{mainthm1} by establishing it for type III groups. In fact, we establish a stronger result which will enable us to complete the proof of Proposition \ref{prop371} (and hence Theorem \ref{mainthm2}) too.
\begin{lemma}\label{lem35} Suppose that $G$ is type III with exponent $m$ and that $A \subseteq G$ has at most $\delta n^2$ Schur triples, where $\delta \leq 2^{-93}m^{-5}$. Then
\[ |A| \; \leq \; \left(\nu(G) + 64m^{2/3}\delta^{1/3}\right)n.\]
\end{lemma}
\proof Recall that $\eta = 2^{-23}$. As $\kappa := 32 \delta^{1/3} m^{2/3} \leq \eta/8m$ and since $q \leq m$, Proposition \ref{thm16} applies. We deduce that either $|A| \leq (\nu(G) - \eta/8)n$, which clearly implies the result, or else that $\alpha_i \leq 2\kappa$ for all $i \notin \{k+1,\dots,5k\}$, the middle of $\mathbb{Z}/q\mathbb{Z}$. Observe that the middle can be partitioned into $2k$ pairs of the form $\{i,2i\}$, where $i$ ranges over $\{k+1,\dots,2k\} \cup \{4k+1,\dots,5k\}$. Using Lemma \ref{prop9pt1}, which tells us that $\alpha_i + \alpha_{2i} \leq 1 + \kappa$, we have $\sum_i \alpha_i \leq 2\kappa q + 2k$, which means (since $\nu(G) \geq 2k/q$), that $|A| \leq n(2\kappa q + 2k)/q \leq (\nu(G) + 64\delta^{1/3}m^{2/3})n$.\endproof
\begin{corollary} $\nu(G) = \mu(G)$ for type III groups.
\end{corollary}
\proof Simply set $\delta = 0$ in the above to get that if $A \subseteq G$ is sum-free then $|A| \leq \nu(G)n$, and so $\mu(G) \leq \nu(G)$. But we have already observed that $\mu(G) \geq \nu(G)$.\endproof
\begin{corollary}\label{prop371b} Proposition \ref{prop371} holds for groups of type III. That is, if $G$ is of type III and if $A \subseteq G$ has $\delta n^2$ Schur triples then $|A| \leq \left(\mu(G) + 2^{20}\delta^{1/5}\right)n$.\end{corollary}
\proof Suppose then that $G$ is of type III and has exponent $m$. Suppose that $A \subseteq G$ has $\delta n^2$ Schur triples. If $\delta \leq 2^{-93}m^{-5}$ then we may apply Lemma \ref{lem35}, and it is easily confirmed that the result holds in this case. If $\delta > 2^{-93}m^{-5}$ then we instead use Corollary \ref{cor44}, obtaining
\begin{eqnarray*}
|A| & \leq & \left(\max\left(\frac{1}{3},\mu(G)\right) + 3\delta^{1/3}\right)n \\ & \leq & \left(\mu(G) + \frac{1}{3m} + 3\delta^{1/3}\right)n \\ & \leq & \left(\mu(G) + 2^{20}\delta^{1/5}\right)n,\end{eqnarray*}
as required.\endproof\\[11pt]
Observe that Lemma \ref{prop371a} and Corollary \ref{prop371b} together imply Proposition \ref{prop371}.\\[11pt]
We are also now in a position to count sum-free sets in groups of type III.
\begin{proposition}\label{047} Theorem \ref{mainthm2} holds for groups of type III. That is, for groups of this type we have $\sigma(G) = \mu(G) + O\left((\log n)^{-1/45}\right)$.
\end{proposition}
\proof Identical to the proof of Proposition \ref{046}, using Corollary \ref{prop371b} in place of Lemma \ref{prop371a}. \endproof\\[11pt]
This result and Proposition \ref{046} together imply Theorem \ref{mainthm2}.
\section{Concluding remarks and open problems.}
\ni The results of this paper (if not all of the methods) are fairly satisfactory, in that the problem of finding $\mu(G)$ is solved for all finite abelian groups, and $|\mbox{SF}(G)|$ is estimated fairly accurately. It would be of some interest to get an asymptotic for this quantity for all abelian groups $G$. This may be very difficult, though it is likely that the methods of \cite{GreenCE} would help in certain cases, particularly $G = C_p$. \\[11pt]
One case of particular interest seems to be $G = C_7^m$, which caused us so much difficulty in the present paper. If $H \leq G$ is any subgroup of index $7$  then we can construct sum-free subsets $A \subseteq G$ as follows. Identify $G/H \cong \mathbb{Z}/7\mathbb{Z}$, let $k = \lfloor \log_2 n\rfloor$ and pick any subset $S \subseteq H + 2$ with $|S| = k$. Let $A$ consist of $S$ together with an arbitrary subset of $H + 3$ and an arbitrary subset of $(H + 4)\setminus (S + S)$. In this way one gets at least $\binom{|H|}{k}2^{2n/7 - k(k+1)/2}$ sum-free subsets of $G$. Using the estimate $\binom{a}{b} \geq (a/eb)^b$, one can easily confirm that this is $\gg 2^{c(\log n)^2}2^{2n/7}$.\\[11pt]
If finding an asymptotic is too ambitious, one could still look to improve on the error term in Theorem \ref{mainthm2}. \\[11pt]
It is possible to generalize the notion of sum-free set to non-abelian groups. If $G$ is non-abelian then we say that $A \subseteq G$ is \textit{product-free} if there are no solutions to $xy = z$ with $x,y,z \in A$. Write $\mu(G)$ for the density of the largest product-free set in $G$. Very little is known concerning $\mu(G)$. Kedlaya \cite{Kedlaya} has shown that $\mu(G) \gg |G|^{-3/14}$, but even the following question is unresolved.
\begin{question}\label{q7.8.1}
Is there a sequence of groups $\{G_n\}_{n = 1}^{\infty}$ with $|G_n| \rightarrow \infty$ and $\mu(G_n) \rightarrow 0$?
\end{question}
\ni For all we know it may be possible to take $G_n = A_n$, the alternating group on $n$ letters. It may even be the case that for any fixed $\epsilon > 0$ there exists $n_0(\epsilon)$ with the following property: if $A_n$ is the alternating group on $n > n_0(\epsilon)$ letters and $A \subseteq A_n$ is a subset of size $|A| \geq \epsilon|A_n|$, then $|AA^{-1}| \geq (1 - \epsilon)|A_n|$.
\section{Acknowledgements.} \ni The authors would like to thank Vsevolod Lev for his extremely careful reading of this paper, which led to a number of improvements in the exposition.

     \end{document}